\documentclass[11pt]{article}
\usepackage{latexsym}  
\usepackage{amsmath}  
\usepackage{amssymb}  
\usepackage{amsfonts}  
\usepackage{verbatim}  
\allowdisplaybreaks[4]  
 \newtheorem{Th}{Theorem}[section]  
 \newtheorem{lem}[Th]{Lemma}  
   
 \newtheorem{pro}[Th]{Proposition}  
   
 \newtheorem{ex}[Th]{Example}  
 \newtheorem{rk}[Th]{Remark}

 \newcommand{\bpr}{\emph{Proof: }}  
\def\epr{ \hfill $\Box$}

\let \ssection=\section  
\renewcommand{\section}{\setcounter{equation}{0}\ssection}  
\newcommand{\e}{\varepsilon}

\newcommand{\bN}{\mathbb{N}}

\newcommand{\bR}{\mathbb{R}}

\title{On strongly Petrovski\u{\i}'s parabolic SPDEs in arbitrary dimension}
\author{C. Cardon-Weber \& A. Millet\\
 \small{ \textit{
Laboratoire de Probabilit\'es et Mod\`eles Al\'eatoires,
  CNRS-UMR 7599}}, \\
\small{ \textit{ Universit\'e Paris 6, 4 place Jussieu,
    Tour 56, F-75252 Paris Cedex 05}}\\
\small{ \textit{ and MODAL'X, Universit\'e Paris 10}}\\
\small{ \textit{ email:  weber@ccr.jussieu.fr and amil@ccr.jussieu.fr}}}
\date{}
 
\begin{document}
 \maketitle

 \begin{abstract}
In this paper we show that the Cahn-Hilliard stochastic SPDE has a
function valued solution in dimension 4 and 5 when the
perturbation is driven by a space-correlated Gaussian noise. This
is done proving general results on SPDEs with globally Lipschitz
coefficients associated with operators on smooth domains of
$\mathbb{R}^d$ which are parabolic in the sense of
Petrovski\u{\i}, and do not necessarily define a semi-group of
operators. We study the regularity of the trajectories of the
solutions and the absolute continuity of the law at some given
time and position.
\end{abstract}

\section{Introduction - Weak solution}  
Let $Q$ be a compact subset of $\bR^d$, %$\mathcal{O}=\stackrel{\circ}{Q}$, 
$\sigma$ and $b_i\, \,  
1\leq i\leq N$ be real-valued functions defined on $[0,T]\times  
Q\times \mathbb{R}$ and  $(k_i, 1\leq  
i\leq N)$ be multi-indices. Let $F$ denote a 
one-dimensional  
(d+1)-parameter Gaussian noise (either a space-time white noise, or a 
space correlated noise),  $A(t,x,D_x)$ denote a differential operator of order $2n$. 
Consider the following  
stochastic partial differential equation  
\begin{equation}\label{weak}  
\frac{\partial}{\partial t} u(t,x)=A(t,x,D_x)\,  
u(t,x)+\sigma(t,x,u(t,x))\, \dot{F}(t,x) + \sum_{i=1}^N D_x^{k_i}  
\big( b_i(t,x,u(t,x))\big)\, ,  
\end{equation}  
with  some homogeneous boundary  conditions, denoted by   \textbf{(BC)}, 
defined for $1\leq q\leq n$ by 
\[
B_q(t,x,D_x)u(t,x)=
\sum_{|k|\leq r_q} B_{q,k}(t,x)\, D_x^k u(t,x)=0\; \mbox{\rm on} \; [0,T]\times \partial Q, \; r_q\leq 2n-1\, ,
\]
and the initial condition $u(0,.)=u_0$. 
We suppose that the operator $L= \frac{\partial}{\partial t} - A(t,x,D_x)$ is uniformly parabolic in the  
sense of Petrovski\u{\i}  and that the boundary conditions are complementary
and satisfy the normality assumption  (for the complete definition see e.g.  
S.D. Eidelman and N.V. Zhitarashu \cite{EZ}, p. 2-17). The most important 
class of uniformly  
strongly parabolic operator in the sense of Petrovski\u{\i} is  
defined by: 
$$A(t,x,D_x)=\sum_{\vert k\vert \leq 2n  
}a_k(x,t)D_x^k\, ,$$ 
where there exists a positive constant  
$\delta_0$ such that for any $(x,t) \in Q\times  
[0,T],\, \xi \in \bR^d$, $$ (-1)^n \Big (\sum_{\vert k\vert =2n}  
a_k(x,t)\:  \xi^k\,  \Big ) \leq -\delta_0\vert  
\xi\vert^{2n}. $$
In this particular case, the definitions can be found in \cite{LM}, p.113-121.
A simple example is provided by $A(t,x,D_x)=a_1(t,x) D_n( a_2(t,x) D_n)$, where $D_n$ is a differential operator of order  $n$, 
and $(-1)^n \sup\{ (a_1a_2)(t,x)\, :\, (t,x)\in [0,T]\times Q\} <0$. When $n=2$ and $D_2=\Delta$, the   Dirichlet boundary conditions ($u=\Delta u=0$ on $[0,T]\times \partial Q$) or Neumann boundary conditions ($\frac{\partial}{\partial \nu}u= \frac{\partial}{\partial \nu} \Delta u=0 $ on $[0,T]\times \partial Q$, where $\nu$ denotes the outer normal)  are normal and complementary. 
\smallskip

A similar equation has  
been studied by Z.~Brzezniak and S.~Peszat \cite{BP} when  
$Q=\mathbb{R}^d$  
 for smooth bounded coefficients  
$a_k$ which do not depend on the time parameter $t$.  
Using time-homogeneous semi-group techniques,  
 these authors prove the existence  
and uniqueness of a mild solution to (\ref{weak})  in some  
weighted $L^p$ spaces (or the set of continuous functions having  
some decay property at infinity), depending on the hypothesis on  
the initial condition $u_0$. Unlike this paper, we allow the case  
where the differential operator does not yield a time-homogeneous  
semi-group and work with martingale-measures as in J.B. Walsh  
\cite{W}.

Let $l\geq 0$ be an integer and $\lambda \in ]0,1[$. According to
S.D.~Eidelman and N.V.~Ivasisen \cite{EI}, if $\partial Q$ is of class
$\mathcal{C}^{2n+l
+\lambda}$, the coefficients $a_k(t,x)$ are of class $\mathcal{C}^{2n(l+\lambda),l+\lambda}([0,T]\times Q)$ for $|k|\leq 2n$, the coefficients $B_{q,k}(t,x)$ are of class $\mathcal{C}^{(2n-r_q+l+\lambda)2n, 
2n-r_q+l+\lambda}([0,T]\times \partial Q)$, then if $G$ denotes the  
 Green  function associated with $L$ and the boundary conditions  
\textbf{(BC)},  for $|a|+2nb\leq 2n+l$,  one has for every $t> s \geq 0$ and
$x,y\in Q$
   \begin{eqnarray}\label{DDH}  
& & \left|D_x^a\, \frac{\partial^b}{\partial t^b}\, G(t,x; s,y)\right|  
 \leq C (t-s)^{-(\alpha +a\, \delta+ b\, \eta)}\,  
 \exp\left(-c\, \frac{|x-y|^{\beta}}{(t-s)^\gamma} \right)\, ,\\  
& & \textrm{with } \alpha =\frac{d}{2n},\, \beta=\frac{2n}{ 2n-1},  
\,\gamma=\frac{1}{2n-1}, \, \delta=\frac{1}{2n} \textrm{ and } \eta=1.\label{DDHbis}  
 \end{eqnarray}  
 Note  
that in some cases, it is possible to extend this upper estimate  
on $G$ to the case the subset $Q$ is not smooth  
(see e.g. \cite{CW1} for the case $A=-\Delta^2$ on $Q=[0,\pi]^d$ and homogeneous Neumann's boundary conditions).  
Therefore, for $d<2n$, the integral $\int_0^t\int_Q G^2(t,x;s,y)\, dy ds$ converges, so that the stochastic integrals of $G(t,x;s,y)$ with respect to the space-time white noise $F(ds,dy)$ are well-defined.
% for any $(t,x)$.
Usual arguments  show that in the particular case 
 of selfadjoint operators $A(t,x,D_x)$, such as $A=a_1(t) D_n( a_2(t,x) D_n)$ with appropriate normal and complementary Dirichlet boundary conditions, the
 Green function $G(t,.;s,.)$ is symmetrical in $x$ and $y$, so that 
 $D_y^aG(t,x;s,y)= D_y^{\tilde{a}}G(t,x;s,y)$ with $|a|=|\tilde{a}|$; then (\ref{DDH}) holds for $D_y^a$ instead of $D_x^a$.\smallskip
 
 We now generalize the setting of \cite{W}, in order to define a "weak" solution
 to (\ref{weak}), which is an alternative to  mild solutions.  Return time, and consider the adjoint operator $L^*=-\frac{\partial }{\partial t} -A^*(t,x,D_x)$ and the adjoint boundary conditions $B^*_q=0$, $1\leq q\leq n$, on $[0,t]\times \partial Q$; then   for fixed $t>0$, exchanging the role of $(t,x)$ and $(s,y)$, $G(t,x;s,y)$ is the fundamental solution to the adjoint problem on the time interval $[0,t]$. Thus, for any smooth function $\phi$ on $Q$, the function 
 \begin{equation}\label{adjoint}
 v(s,y)=\int_Q G(t,x;s,y)\, \phi(x)\, dx
 \end{equation}
  is the solution to the equation $L^*v=0$ on $[0,t]\times Q$, with
  adjoint boundary conditions $(B^*_q v=0, 1\leq q\leq n)$ on
  $[0,t]\times \partial Q$, and such that $v(t,.)=\phi$. Then for
  Dirichlet's systems ($r_q=q-1$, $1 \leq q \leq n$) or in particular
  cases (see e.g. example \ref{Excond}), for $v$ defined by
  (\ref{adjoint}) and ``regular'' $u$, the following Green formula holds 
   (see e.g. \cite{EZ}, p. 231 or \cite{LM}, p. 133):
 \begin{equation}\nonumber
 \int_0^t \int_Q (Lu)(s,y)\, v(s,y)\, dy ds +\int_Q u(0,y)\, v(0,y)\, dy =
 \label{formgreen}   \int_Q u(t,x)\, v(t,x)\, dx \, .
 \end{equation}
  Furthermore, if
 integration by parts yields 
 \[ \int_0^t\int_Q D_y^{k_i}(b_i(s,y,u(s,y))\, v(s,y))\, dy ds = (-1)^{|k_i|} \int_0^t\int_Q b_i(s,y,v(s,y))\, D_y^{k_i}v(s,y)\, dy ds \]
  and if  $\int_0^t\int_Q \sigma(s,y,u(s,y))\, v(s,y) F(ds,dy)$ is defined as the stochastic integral with respect to a worthy martingale measure, then even if $u$ is not "regular", we can define the weak solution to (\ref{weak}),  by requiring that the following form of equation (\ref{formgreen}) holds for every function $v$ in $\mathcal{C}^{2n,1}([0,T]\times Q)$:
 \begin{eqnarray}\nonumber 
&& \int_0^t \int_Q \left[\sigma(s,y,u(s,y))\, v(s,y)\, F(ds, dy) +
 \sum_{i=1}^N (-1)^{|k_i|} b_i(s,y,u(s,y)) D_y^{k_i}v(s,y)\, dy ds \right]\\
 &&\qquad \qquad +\int_Q u(0,y)\, v(0,y)\, dy  =  \int_Q u(t,x)\, v(t,x)\, dx\, \label{greenintegre} \, .
\end{eqnarray} 
Using  (stochastic) Fubini's theorem we obtain the following evolution equation, which is equivalent to 
 (\ref{greenintegre}):
 \begin{eqnarray}\nonumber  
u(t,x)&=&\int_Q G(t,x;0,y)\,  u_0(y)\, dy + \int_0^t\int_Q G(t,x;s,y)\,  
 \sigma(s,y,u(s,y))\, F(ds,dy)\\  
&&\qquad + \sum_{i=1}^N \int_0^t\int_Q H_i(t,x;s,y)\, b_i(s,y,u(s,y))\, dy \, ds \, ,  
\label{evol}  
\end{eqnarray}  
where $H_i(t,x;s,y)= (-1)^{|k_i|}D_y^{k_i}G(t,x;s,y)$; if $G(t,x;s,y)$ is
symmetric in $x$ and $y$, the upper estimate  
 (\ref{DDH}) implies that  
%\[
$|H_i(t,x;s,y)|\leq C\,  (t-s)^{-(\alpha + |k_i|\delta) }\,  
 \exp\left(-c\, \frac{|x-y|^{\beta}}{(t-s)^\gamma} \right).$
 %\]  
 We give an example where all the requirements (except that on the existence
 of the stochastic integral) are fulfilled.
 \begin{ex}\label{Excond}
 The boundary of the set $Q$ is of class $\mathcal{C}^{4+l+\lambda}$,
 the functions $a_1(t)\in\mathcal{C}^{4(l+\lambda)}([0,T])$, $a_2(t,x)\in \mathcal{C}^{4(l+\lambda),l+\lambda}([0,T]\times Q)$, $A(t,x,D_x)= a_1(t)\, \Delta (a_2(t,x)\, \Delta )$ and ${\displaystyle \sup_{(t,x)\in [0,T]\times Q} a_1(t)\, a_2(t,x) <0}$.
 
 {\bf Case 1} Dirichlet's boundary conditions: $u=\Delta u=0$ on $[0,T]\times \partial Q$,
either  $0\leq |k_i|\leq 1$, or $|k_i|=2$ and $b_i(s,y,z)=\tilde{b}_i(z)$ for some
 function $\tilde{b}_i$ of class $\mathcal{C}_2$ such that $\tilde{b}_i(0)=0$.

{\bf Case 2} Neumann's boundary conditions: let $\nu$ denote the outer normal, 
$\frac{\partial}{\partial \nu}u= \frac{\partial}{\partial \nu}\Delta u=0$, 
$\frac{\partial}{\partial \nu}a_2(t,x)=0$ on $[0,T]\times \partial Q$, $k_i=0$ or
$|k_i|=2$, $k_i$ has even components and $b_i(s,y,z)=\tilde{b}_i(z)$ for some
 function $\tilde{b}_i$ of class $\mathcal{C}_2$.
\end{ex}

For $d\geq 2n$, the function $G^2(s,x;s,y)$ need  
not be in $L^2([0,T]\times Q, dsdy)$, so that the Gaussian noise  
$F$ need not be the space-time white noise; we require the noise $F$ to be a Gaussian process which is  
white in time, but has a space correlation defined in terms of a  
function $f$ depending on the difference of two vectors of  
$\mathbb{R}^d$ (or such that when $|x-y|\rightarrow 0$, the product
$f(x,y)\, |x-y|^a$ remains bounded for some $a>0$). 
 We just mention a few previous papers on this  
subject, stressing the type of noise which is used. A particular  
case of this noise (where the function $f$ only depends of the  
norm  $|x-y|$, such as $f(x-y)=|x-y|^{-a}$ for $0<a<d$)  has  
been used in  C.~Mueller \cite{Mu}, R.~Dalang and N.~Frangos  
\cite{DF}, A.~Millet and M.~Sanz-Sol\'e \cite{MS} in the case  
$Q=\mathbb{R}^2$ for the  wave operator. In these papers, the  
existence and uniqueness of a continuous solution is proved by  
precise estimates of integrals involving the corresponding Green  
function.  A more general covariance structure (depending on the  
radon measure $\mu$ with Fourier transform $f$, or more generally
a tempered distribution $\Gamma=\hat{\mu}$) has been used in  
\cite{BP}, \cite{KZ1}, \cite{KZ2}, \cite{PZ2},\cite{P}  for the wave 
and heat  
operators on $\mathbb{R}^d$; in the last references, the existence  
of a solution is proved in some weighted $L^p$-space, or in the  
space of continuous functions with some decay at infinity, and the  
method uses the Fourier transform of the Green function  
$G$ (see e.g. S.~Peszat and J.~Zabczyk \cite{PZ2} for a detailed account of existence
and uniqueness results to parabolic SPDEs with a semi-group structure
in any dimension).  
  This general covariance was also used by R.~Dalang \cite{D}, who proves  
  the existence of continuous processes solutions to the heat and wave  
  stochastic SPDEs by means of an extension of stochastic  
 integrals with respect to martingale measures for distribution-valued  
 integrands.  In these references, the coefficients of the
differential operator $A$
do not depend on $(t,x)$.
\medskip

On the other hand, several attempts have been made to find function-valued  
solutions to  "highly non-linear"  
stochastic SPDEs, namely PDEs  with a polynomial forcing term  $b_i$  
 (such as the Burgers  PDE ($d=1$, $A=\Delta$, $N=1$, $a_1=1$ and $b_1(t,x,y)=y^2$),  
 or  the  Cahn-Hilliard's PDE ($d\leq 3$, $A=-\Delta^2$, and  
 $\sum_i D_x^{k_i}b_i(t,x,u(t,x)) =  
\Delta R(u(t,x))$, where $R$ is a polynomial of odd degree with  
positive dominant coefficient)  and with  
 a stochastic perturbation driven by the space-time white noise.  
Thus, G.~Da Prato A.~Debussche and Temam \cite{DPDT} and then I.~Gy\"ongy  
\cite{Gy} have proved the existence of a  
 function-valued solution to the stochastic Burgers  equation in dimension  
 1.  G.~Da Prato and A.~Debussche \cite{DPD} have proved the existence  
of a function-valued solution to the stochastic Cahn-Hilliard  
equation in dimension 1 (up to 3) when the perturbation is driven  
by a space-time white noise (a Gaussian noise with some spatial  
correlation).  
  C.~Cardon-Weber \cite{CW1} and  
 \cite{CW2} has proved the existence of a function-valued solution to the stochastic  
 Cahn-Hilliard equation in dimension $d\leq 3$  when $R$ is a  
 polynomial of degree 3  and when  the stochastic  perturbation is driven by the  
 space-time white noise.  The method  
 used in these papers is the following: using a truncation procedure and the  
 existence and uniqueness  
results proved in the case of globally Lipschitz coefficients, one  
proves the existence and uniqueness of a solution to the SPDE  
where the polynomial coefficients have been changed.  
 Then the uniqueness property of the solution allows to use concatenation to  
  obtain the existence of  a solution up to some stopping time. Finally,  a  
  priori estimates for a deterministic PDE obtained by isolating the  
  stochastic integral (whose behavior is  
 controlled by means of the Garsia lemma), prove that this stopping time is  
 the terminal time $T$. These last estimates use methods of analysis which  
heavily depend on the specific form of the PDE, and no general  
scheme can be given. Let us finally mention that, using semi-group  
techniques, Z.~Brzezniak and S.~Peszat \cite{BP} have proved the  
existence of solutions to some SPDE with a polynomial drift term  
(when $N=1$, $a_1=0$, and when the operator $A$ is of order 2 and  
yields a semi-group of operators). Also in \cite{MM}, the  
existence of the solution to a stochastic wave equation in  
dimension 2 with a non-uniformly  Lipschitz drift has been proved,  
while the existence to the stochastic KDV equation has been shown  
by A.~Debussche and A.~de Bouard \cite{DDB}.  
\bigskip

The aim of this paper is two-fold. On one hand, we prove that in this
 general context with time and space dependent coefficients,  the  
upper estimates (\ref{DDH})  
  of the Green function $G$ and its time and space derivatives are sufficient  
   to ensure  
the existence and uniqueness of the solution $u$ to (\ref{evol}),  
provided that some integrability condition of the covariance  
function $f$ on a neighborhood of 0 is required.  
 We prove that, when the Green function $G$  has a lower estimate by  
 $t^{-\alpha}$ on the  diagonal  (which can be  the case when it admits  
 an explicit eigenvectors-eigenvalues expansion), this condition is necessary  
to be able to consider stochastic integrals of $G$. As in \cite{W}  
and \cite{D}, we use stochastic integrals with respect to  
martingale measures. We give sufficient conditions on  
 the covariance function $f$ for  the trajectories of $u$ to be  
 H\"older-continuous. We then use these  results  to extend the existence  
 and uniqueness of a  
function-valued solution to the stochastic Cahn-Hilliard equation  
when $Q=[0,\pi ]^d$ or a bounded "smooth" subset of $\mathbb{R}^d$, $d=4, 5$. We give necessary and  
sufficient conditions on the covariance  
 function $f$ to ensure that the stochastic integral of the corresponding  
 Green function $G$ is well-defined. We study regularity properties of the  
 trajectories of $u$ and prove that,  
if $Q=[0,\pi]^d$ and the diffusion coefficient $\sigma$ is strictly elliptic,  the  
law of $u(t,x)$ has a density for $t>0$ and $x\in Q$. This  
extends the results proved in \cite{DPD} and \cite{CW1} to higher  
dimensions.  For the sake of simplicity, we mostly restrict ourselves to the
case $F$ is a space-correlated noise; this could be avoided in "small" dimension
for arbitrary Petrovski\u{\i}'s parabolic SPDEs. Also note that the proof of the 
existence of a solution to the stochastic Cahn-Hilliard equation extends 
directly to the more general situation described in Example \ref{Excond},
when $R$ is a polynomial of degree 3 with  positive dominant coefficient (which has
no constant term in the case of the Dirichlet boundary conditions).
 \bigskip

The paper is organized as follows. Section 2 gives necessary and sufficient  
conditions  
to ensure that stochastic integrals of a function $G$ which satisfies  
(\ref{DDH}) are well  
defined, and sufficient conditions to ensure H\"older properties of  
stochastic integrals  
appearing in (\ref{evol}), provided that the process $u$ has bounded moments.  
 In section 3, we prove both the existence of solutions to (\ref{evol})  
 either in $\mathcal{C}([0,T], L^q(Q))$  
 or $\mathcal{C}([0,T]\times Q)$ when the coefficients are globally  
 Lipschitz functions.  
 We then concentrate on the proof of the a priori estimates which allow to  
 deduce the existence and uniqueness of a solution to the stochastic  
Cahn-Hilliard equation in dimension 4 and 5.  
  Section 4  establishes  H\"older regularity of the trajectories, while  
  section 5 shows the absolute continuity of the law of $u(t,x)$ for  
  $t>0$ and $x\in Q$ and the solution $u$ to the Cahn-Hilliard equation.  
  \bigskip

All the constants $C$ appearing in the statements can change from  
one line to the next one. When we want to stress the fact that  
$C$ depends on some parameter $k$, we denote it by $C_k$.

\section{Stochastic integrals with respect to a space correlated noise}  
Let $Q$ be a compact subset of $\mathbb{R}^d$,  
$\mathcal{D}(\mathbb{R_+}\times Q)$ denote the space of functions  
$\varphi\in \mathcal{C}_0^\infty( \mathbb{R}_+\times Q)$ with  
compact support,  endowed with the topology defined by the  
following  convergence: $\varphi_n \rightarrow \varphi$ if:  
\medskip

(i) There exists a compact subset $K$ of $\mathbb{R}^+\times Q$ such that ${\rm support}  
(\varphi_n-\varphi) \subset K$ for all $n$.  
  
(ii) $\lim_{n\rightarrow +\infty} D^a \varphi_n = D^a \varphi$ uniformly  
on $K$ for every multi-index $a$.  
\medskip

Let $F=\big(F(\varphi)\, , \, \varphi\in  
\mathcal{D}(\mathbb{R}_+\times Q)\big)$ be an $L^2(P)$-valued  
centered Gaussian process, which is white in space but has a space  
correlation defined as follows: given $\varphi$ and $\psi$ in  
$\mathcal{D}(\mathbb{R}_+\times Q)$,  
 the covariance functional of $F(\varphi)$ and $F(\psi)$ is  
\begin{equation}\label{cov}  
J(\varphi,\psi)=E\big( F(\varphi)\, F(\psi)\big) =  
\int_0^{+\infty}dt\int_Qdy\int_Q\varphi(t,y)\, f(y-z)\, \psi(t,z)\, dz\,  
\end{equation}  
where  $(Q-Q)^*=\{y-z\, : y,z\in Q, \, y\neq z\}$ and  
$f:(Q-Q)^*\rightarrow [0,+\infty[$ is a continuous function.  
According to \cite{S}, the bilinear form $J$ defined by  
(\ref{cov})is non-negative definite if and only if $f$ is the  
Fourier transform of a non-negative tempered distribution $\mu$ on  
$Q$. Then $F$ defines a martingale-measure (still denoted by $F$),  
which allows to define stochastic integrals (see \cite{W}).  
  
In this section, we consider fairly general functions  
$H:([0,T]\times Q)^2\rightarrow \mathbb{R}$, including Green  
functions associated with parabolic operators in the sense of  
Petrovski\u{\i}: more precisely, we suppose that $H$ satisfies the  
following upper estimate for  
 some  $\beta\geq 1$,   some strictly positive parameters  
 $\alpha $ and $\gamma$ and some  positive constants $c$ and $C$: for  
  any $t>0$, $x,y\in Q$:  
\begin{equation}  
\label{majoH}  
|H(t,x;s,y)|\leq C\, (t-s)^{-\alpha}\, \exp\left( -c\frac{|x-y|^\beta}  
{(t-s)^\gamma}\right)\, .  
\end{equation}  
Then the change of variables $u=(y-x)\, (t-s)^{-\frac{\gamma}{\beta}}$  
yields  
\begin{equation}\label{intDH}  
\int_Q|H(t,x;s,y)|\, dy\leq C\, (t-s)^{-\alpha +\frac{\gamma}{\beta}d}\,  
.  
\end{equation}  
  
We now give a sufficient integrability condition on the  
space-correlation function $f$ of the Gaussian noise $F$ to ensure  
that the stochastic integral of a  bounded adapted process  
multiplied by  a kernel satisfying (\ref{majoH})  is a  
well-defined stochastic process. Under some additional assumptions  
on $H$ we prove that this condition is necessary.  
\begin{lem}.\label{CNS}
(i) Let $H :([0,T]\times Q)^2\rightarrow \bR$ satisfy (\ref{majoH}) and suppose that either  
\begin{equation}\label{Acovd=}  
\int_{B_d(0,1)} f(v)\, \, \ln\left(|v|^{-1}\right)\, dv <\infty,  
\quad if \quad d=\frac{\beta}{\gamma}(2\alpha -1)\, ,  
\end{equation}  
or  
\begin{equation}\label{Acovdiff}  
\int_{B_d(0,1)} f(v)\, |v|^{-[\frac{\beta}{\gamma}(2\alpha -1)-d]^+}\,  
 dv <\infty, \quad if \quad  d\neq \frac{\beta}{\gamma}(2\alpha -1)\, .  
\end{equation}  
Then  for any $t\in [0,T]$, 
\begin{equation}  
\label{intconv} I(t)=\int_0^T \int_Q\int_Q |H(t,x;s,y)|\, f(y-z)\,  
|H(t,x;s,z)|\, dy\, dz\, dt <+\infty\, .  
\end{equation}  
(ii) Let $H :([0,T]\times Q)^2\rightarrow \bR$ satisfy  
\begin{equation}\label{minodiag}  
\inf\{H(t,x;s,x)\, :\, x\in Q\}\geq C_0 (t-s)^{-\alpha}  
\end{equation}  
and for every multi-index $k$ with $|k|=1$:  
\begin{equation}\label{majDyG}  
\sup\big\{ \big| D^k_y H(t,x;s,y)\big|\, :\, x,y\in Q\big\} \leq C_1  
(t-s)^{-(\alpha +\delta)}\, ,  
\end{equation}  
for $\alpha, \delta >0$, some constant $C>0$ and any $t>s>0$.  Then  
if (\ref{intconv}) holds for every $x\in Q$ and  $t\in ]0,T]$,  one has either  
(\ref{covdelta=}) or (\ref{covdeltadif}) depending on $d$, where:  
\begin{equation}\label{covdelta=}  
\int_{B_d(0,1)} f(v)\, \, \ln\left(|v|^{-1}\right)\, dv <\infty,  
\quad if \quad d=\frac{1}{\delta}(2\alpha -1)\, ,  
\end{equation}  
or  
\begin{equation}\label{covdeltadif}  
\int_{B_d(0,1)} f(v)\, |v|^{-[\frac{1}{\delta}(2\alpha -1)-d]^+}\,  dv  
 <\infty, \quad if \quad  d\neq \frac{1}{\delta}(2\alpha -1)\, .  
\end{equation}  
\end{lem}  
\begin{rk}.\label{petrov} When $H$ is the Green function of the
operator $\frac{\delta}{\delta t}+\Delta^2$ on $Q=[0,\pi]^d$ 
with homogeneous Neumann's or Dirichlet's boundary conditions, condition
(\ref{minodiag}) holds with $\alpha=\frac{d}{4}$. We only sketch the proof
for Neumann's boundary conditions; in that case, 
$G(t,x;s;x)\geq \sum_{k\in \mathbb{N}^{*d}} \exp(-|k|^4 (t-s))\left( 
\prod_{1\leq i\leq d}\cos^2(k_i x_i)\right)$. Set
$\mathcal{I}= ]\frac{\pi}{2}-\theta, \frac{\pi}{2}+\theta[$ for
some $\theta>0$ small enough; then for $k\in \mathbb{N}^*$
and $x\in ]0,\pi[$, $kx \in \mathcal{I}$ implies that 
$(k+1)x\not\in \mathcal{I}$. Skipping at most every other term
and using the monotonicity of $(\exp(-k^4 (t-s)); k\geq 1)$,
we deduce that $G(t,x;s,x)\geq C \sum_{k\in \mathbb{N}^{*d}} 
\exp(-c|k|^4 (t-s))= C (t-s)^{-\frac{d}{4}}$ for some positive constants
$c$ and $C$ (see Remark 3.5 in \cite{CW1}
for a similar argument).
Since in the case of Green functions of parabolic operators, 
$\delta=\frac{\gamma}{\beta}$,  the conditions  
(\ref{Acovd=}) and (\ref{covdelta=}) (respectively  
(\ref{Acovdiff}) and (\ref{covdeltadif})) are identical.  
\end{rk}  
{\it Proof of Lemma \ref{CNS}:} (i) Since $\beta \geq 1$, $|y-z|^\beta \leq 2^{\beta  
-1}(|x-y|^\beta +|x-z|^\beta)$, so that for $\lambda\in ]0,1[$ and $t_0\in [0,T]$,
\begin{eqnarray}I(t_0)&\leq& C \int_0^T  \, t^{-2\alpha} dt \int_Q  
\exp\left(-c\, (1-\lambda)\frac{|x-y|^\beta}{t^\gamma}\right) dy  
 \nonumber \\  
&&\qquad\qquad\times \int_Q  
\exp\left(-c\, 2^{1-\beta}\, \lambda\frac{|y-z|^\beta}{t^\gamma}\right) \,  
f(y-z)\,dz \nonumber \\  
&\leq& C\, \int_0^T  \, t^{-2\alpha} dt \int_Q  
 \exp\left(-\bar{c}\frac{|x-y|^\beta}{t^\gamma}\right) dy \int_Q  
\exp\left(-\bar{c}\,  \frac{|y-z|^\beta}{t^\gamma}\right) \,  
f(y-z)\,dz \;.\qquad\label{numero}  
\end{eqnarray}  
Set $t^{\frac{\gamma}{\beta}}\eta = x-y$, $v=y-z$ and then $u=|\eta|$;  
 there exist positive constants $C$, $c$ and $R$ such that  
\begin{eqnarray}  
I(t_0)&\leq& C\int_0^T  \, t^{\frac{\gamma}{\beta}d-2\alpha} dt  
\int_0^{+\infty} \exp\left( -c u^\beta\right) \, u^{d-1} \,  
du\int_{B_d(0,R)} \exp\left( -c\,  
\frac{|v|^\beta}{t^\gamma}\right) \, f(v)\, dv\nonumber\\ &\leq&  
C\int_0^T \psi(t)\, dt \, , \label{inter}  
\end{eqnarray}  
 where for any $0<t\leq T$ one sets  
\begin{equation}  
\label{defpsi} \psi(t)=t^{\frac{\gamma}{\beta}d-2\alpha}  
\int_{B_d(0,R)} \exp\left( -c\frac{|v|^\beta}{t^\gamma}\right)\,  
 f(v)\, dv\, .  
\end{equation}  
For fixed $v\neq 0$, set $r=|v|^\beta\, t^{-\gamma}$; then  
Fubini's theorem yields  
\[ \int_0^T \psi(t)\, dt \leq C\int_{B_d(0,R)} f(v)\,  
|v|^{d+\frac{\beta}{\gamma}(1-2\alpha)}\, dv \int_{|v|^\beta\,  
T^{-\gamma}}^{+\infty} \, r^{-1+\frac{1}{\gamma}(2\alpha -1  
-\frac{\gamma d}{\beta})}\, \exp(-c r)\, dr\, .\] We now  
distinguish three cases:  
\smallskip

{\em Case 1} If $2\alpha > 1+ \frac{\gamma d}{\beta}$, then the  
second integral is bounded by a constant independent of $|v|$ and  
$\int_0^T \psi(t)\, dt \leq C \int_{B_d(0,R)} f(v) \,  
|v|^{d+\frac{\beta}{\gamma}(1-2\alpha)}\, dv\, ,$ which yields  
(\ref{Acovdiff}).  
\smallskip

{\em Case 2} If $2\alpha =1+ \frac{\gamma d}{\beta}$, since for  
$0< |v|\leq R$,  
\[\int_{|v|^\beta\, T^{-\gamma}}^{+\infty}  r^{-1}\, \exp(-cr)\, dr  
\leq C\, \vert 1-\ln(|v|)1_{\{\vert v\vert \leq 1\}}\vert   ,\]
 we have $ \int_0^T\psi(t)\, dt \leq C\int_{B_d(0,R)}f(v) \,  
\big[1+ \ln\big(|v|^{-1}\big) 1_{\{\vert v\vert \leq 1\}}\big]\, dv\, ,$ which yields  
(\ref{Acovd=})  
\smallskip

{\em Case 3} Finally, if $2\alpha < 1+\frac{\gamma d}{\beta}$,  
then  for $0< |v|\leq R$,  
\[ \int_{|v|^\beta\, T^{-\gamma}}^{+\infty}  
r^{-1+\frac{1}{\gamma}(2\alpha-1-\frac{\gamma d}{\beta})}\,  
\exp(-cr)\, dr\leq C \left(1+|v|^{\frac{\beta}{\gamma}(2\alpha  
-1)-d}\right)\, ,\] so that $\int_0^T\psi(t)\, dt \leq  
C\int_{B_d(0,R)}f(v) \,  dv\, ,$ which yields (\ref{Acovdiff}).  
\smallskip

Note that for small $T$, the following computation gives a more  
precise upper estimate of $\int_0^T \psi(t)\, dt$, which will be  
used in the sequel. Indeed, for $\nu\in ]0,\gamma[$, the  
decomposition of the integral over $B_d(0,R)$ into $\{ |v|^\beta  
T^{-\gamma} \geq T^{-\nu}\}$ and its complement yields  
\begin{eqnarray*}  
\int_0^T \psi(t)\, dt& \leq &C \int_{B_d\left(0, T^{\frac{\gamma  
-\nu}{\beta}}\right)} f(v)\,  
|v|^{d+\frac{\beta}{\gamma}(1-2\alpha)}\, dv \int_{|v|^\beta  
T^{-\gamma}}^{+\infty} r^{-1+\frac{1}{\gamma}(2\alpha  
-1-\frac{\gamma d}{\beta})}\, \exp(-c r)\, dr\\  
&& +C\, I\,  
\exp(-\bar{c}\, T^{-\nu} )\, .\end{eqnarray*}  Thus for  
$0<\nu<\gamma$ and $ 2\alpha\neq 1+\frac{\gamma d}{\beta}$ one has  
\begin{equation}  
\label{intpsiTdiff} \int_0^T \psi(t)\, dt \leq C\, \left[  
\exp(-\bar{c}\, T^{-\nu}) + \int_{B_d\big(0, T^{\frac{\gamma  
-\nu}{\beta}}\big)} f(v)\, |v|^{\big[  
-\frac{\beta}{\gamma}(2\alpha -1)-d\big]^+}\, dv \right]\, ,  
\end{equation} while for $0<\nu<\gamma$ and $ 2\alpha = 1+  
\frac{\gamma d}{\beta}$ one has  
\begin{equation} \label{intpsiT=} \int_0^T \psi(t)\, dt \leq C\,  
\left[ \exp(-\bar{c}\, T^{-\nu}) + \int_{B_d\big(0,  
T^{\frac{\gamma -\nu}{\beta}}\big)} f(v)\, \ln\big(|v|^{-1}\big)  
\, dv \right]\, .  
\end{equation}  
\smallskip

(ii) The assumptions (\ref{minodiag}) and (\ref{majDyG}) imply  
that for $|x-y|\leq 2\, C_2\, t^\delta$ with $C_2<\frac{C_0}{4C_1}$  small  
enough, one has  
\begin{eqnarray*}  
H(t,x;s,y)&\geq&H(t,x;s,x)-|H(t,x;s,x)-H(t,x;s,y)|\\  
&\geq& (C_0-2\,C_1\,C_2)\, (t-s)^{-\alpha}\geq \frac{C_0}{2}\, (t-s)^{-\alpha}\, .  
\end{eqnarray*}  
Let $a>0$ be such that $Q_a=\{x\in Q\, :\, d(x,\partial Q) >a\}  
\neq \emptyset$; then for $0<2\, C_2\, t_0^\delta \leq a$,  
$0<s\leq t_0\leq T$, $x\in Q_{2\, C_2\, t_0^\delta}$, $y\in B_d(x,C_2\,  
s^\delta)$ and $z\in B_d(y,C_2\, s^\delta)$, one has $y,z\in Q$.  
Thus Fubini's theorem implies for $x\in Q_{2\, C_2\,  
t_0^\delta}\neq \emptyset$:  
\begin{eqnarray*}  
I(t_0)&\geq&\int_0^{t_0} ds\int_{Q\cap B_d(x,C_2s^\delta)} dy  
\int_{Q\cap B_d(y,C_2s^\delta)}  |H(t_0,x;s,y)|\, f(y-z)\, |H(t_0,x;s,z)|  
dz \\ &\geq&C\, \int_0^{t_0} \, s^{-2\alpha+ d \delta}\, ds  
\int_{B_d(0,C_2s^\delta)} f(v)\, dv\geq C\, \int_{B_d(0, R)} f(v)  
\, dv\int_{\left(\frac{|v|}{C_2}\right)^{\frac{1}{\delta}}}^{t_0}  
s^{-2\alpha +d \delta } ds  
\end{eqnarray*}  
for $R=C_2\, t_0^\delta$. Again we have to study three cases  
depending on the power of $s$.  
\smallskip

{\em Case 1} If $d \delta +1< 2\alpha $, let $\bar{R}=\frac{R}{2}$;  
 then one has  
\begin{equation}\label{minoCNSdiff} I(t_0)\geq C\, \int_{B_d(0, \bar{R})}  
  f(v) \,dv  
\int_{\left(\frac{|v|}{C_2}\right)^{\frac{1}{\delta}}}^{\left(  
\frac{2|v|}{C_2}\right)^{\frac{1}{\delta}}}s^{-2\alpha+d\delta}\,  
ds \geq C\, \int_{B_d(0, \bar{R})}  
f(v)\,|v|^{\frac{1-2\alpha}{\delta}+d}\, dv\, ,  
\end{equation}  
 which yields (\ref{covdeltadif}).  
\smallskip

{\em Case 2} If $d \delta +1= 2\alpha $, let $\nu >0$ and let  
$\bar{R}=R\wedge 1\wedge (C_2^{-\nu}R^{1+\nu}))$; then $|v|\leq  
 \bar{R}$ and  $t_0\leq C_2^{\frac{1}{\delta}}$  imply  
$\left( \frac{|v|}{C_2}\right)^{\frac{1}{(1+\nu)\delta}}\leq t_0$  
and $|v|^{\frac{1}{ 2}}\, C_2^{-1}\leq 1$; hence $C_2\, |v|^{-1}\geq 
|v|^{-\frac{1}{2}}$ and  
\begin{equation}\label{minoCNS=} I(t_0)\geq C\, \int_{B_d(0,\bar{R})}  
f(v) \, dv  
\int_{\left(\frac{|v|}{C_2}\right)^{\frac{1}{\delta}}}^{  
 \left(\frac{|v|}{C_2}\right)^{\frac{1}{(1+\nu) \delta}} } s^{-1}\, ds \geq  
 C\, \int_{B_d(0,  
\bar{R})} f(v) \, \ln\big(|v|^{-1}\big) dv\, ,  
\end{equation}  
 which yields (\ref{covdelta=}).  
\smallskip

{\em Case 3} Finally, if $d \delta +1> 2\alpha $, one has $I(t_0)\geq  
C\, \int_{B_d(0, R)} f(v) \ \, dv\, ,$ which yields  
(\ref{covdeltadif}). \epr  
\bigskip

The following lemma gives sufficient conditions on the covariance  
function $f$ to obtain moment estimates of stochastic integrals  
which yield H\"older regularity of the corresponding process.  For  
this, we impose   an upper estimate of the space and time partial  
derivatives of the kernel $H$: there exist positive constants  
$\delta, \eta, c, C$ such that for any $t>0$, $x,y\in Q$ and $k\in  
\mathbb{N}^d$ with $|k|=1$,  
\begin{equation}  
\label{majoDH} \left|D_x^k H(t,x;s,y)\right| \leq C (t-s)^{-(\alpha  
+\delta)}\, ,  
 %\exp\left(-c\, \frac{|x-y|^{\beta}}{t^\gamma} \right)\, ,  
 %\end{equation}  
%\begin{equation}  
\quad \left|\frac{\partial}{\partial t}H(t,x;s,y)\right| \leq C  
(t-s)^{-(\alpha +\eta)}\,  
 \exp\left(-c\, \frac{|x-y|^{\beta}}{(t-s)^\gamma} \right)\, .  
\end{equation}  
In order to deal with time increments, we  impose also that  
$f$ satisfies the following "monotonicity" condition:  
  
\noindent \textbf{(C1)} There exist strictly  positive constants  
$C_1$ and $c_1$ such that  
     \begin{equation}\label{fmonoton}  
     f(u)\leq C_1 f(v)\quad {\rm for} \quad |v|\leq c_1\, |u|\, .  
     \end{equation}  
Note that (C1) holds if $f(u)=|u|^{-a}$ for some $a>0$.  
\begin{lem}.\label{ISaccrois}  
Suppose that $Q$ is convex and  let  $H : ([0,T]\times Q)^2\rightarrow \mathbb{R}$ satisfy the condition  
(\ref{majoH}). Let $F$ be  
a Gaussian noise with spatial covariance defined by (\ref{cov})  
such that the  
 correlation function $f$ satisfies {\textbf (C1)} and either (\ref{Acovd=}) or 
 (\ref{Acovdiff}).  
Fix $p\in [1,+\infty[$, let  $u:\Omega\rightarrow \mathbb{R}$ be an  
adapted process such that $\displaystyle \sup_{(t,x)\in [0,T]\times Q}  
 E(|u(t,x)|^{2p})<+\infty$  
 and  
 for $t\in[0,T]$ and $x\in Q$, let  
\begin{equation}\label{ISu}  
I(t,x)=\int_0^t \int_Q H(t,x; s,y)\, u(s,y) \, F(ds,dy)\, .  
\end{equation}  
  
(i) Suppose  that $H$ satisfies (\ref{majoDH}) and let $a\in]0,1[$;  
if either  
\begin{equation}\label{covISxd=}  
\int_{B_d(0,1)} f(v)\, \ln\left(|v|^{-1}\right)\, dv <+\infty  
\quad {\rm for}\quad d=\frac{\beta}{\gamma}\,  
(2\alpha+a\delta-1)\, ,  
\end{equation}  
or  
\begin{equation}\label{covISxdif}  
\int_{B_d(0,1)} f(v)\,  
|v|^{-\left[\frac{\beta}{\gamma}(2\alpha+a\delta-1)-d\right]^+} \,  
dv <+\infty \quad {\rm for}\quad d\neq\frac{\beta}{\gamma}\,  
(2\alpha+a\delta-1)\, ;  
\end{equation}  
then  there exists $C_p>0$  such that for  
every $x,x'\in Q$,  
\begin{equation}\label{acroiISx}  
A(x,x')= \sup_{t\in [0,T]} E\left(|I(t,x)-I(t,x')|^{2p}\right)\leq  
C_p\, |x-x'|^{ap}\, .  
\end{equation}  
  
(ii) Suppose that  $H$  
satisfies (\ref{DDH})and (\ref{DDHbis}) and let $b\in]0,1[$ ;  if either  
\begin{equation}\label{covIStd=}  
\int_{B_d(0,1)} f(v)\, \ln\big(|v|^{-1}\big)\, dv <+\infty \quad  
{\rm for}\quad d=\frac{\beta}{\gamma}\, (2\alpha+b\eta -1)\, ,  
\end{equation}  
or  
\begin{equation}\label{covIStdif}  
\int_{B_d(0,1)} f(v)\, |v|^{-[\frac{\beta}{\gamma}(2\alpha+b\eta -1)-d]^+}\,  
 dv <+\infty \quad {\rm for}\quad d\neq\frac{\beta}{\gamma}\,  
 (2\alpha+b\eta -1)\, ;  
\end{equation}  
then  there exists  $C_p>0$ such that for every $0\leq t<t'\leq  
T$,  
%\begin{eqnarray}\label{acroiISt}  
\begin{equation}\label{acroiISt}  
\nonumber B(t,t')= \sup_{x\in Q} E\left(\left|\int_0^t \int_Q  
\Big[H(t,x; s,y)-H(t',x;s,y)\Big]\, u(s,y)\, F(ds,dy)  \right|^{2p}  
\right)%\\ &\leq &  
\leq C_p\, |t-t'|^{bp}\, ,  
\end{equation}  
%\end{eqnarray}  
and  
\begin{equation}\label{IStt'}  
C(t,t')= \sup_{x\in Q} E\left(\left|\int_t^{t'} \int_Q  
H(t',x;s,y)\, u(s,y)\, F(ds,dy) \right|^{2p} \right)\leq C_p\,  
|t-t'|^{bp}\, .  
\end{equation}  
\end{lem}  
\bpr (i) Burkholder's inequality implies that for every $p\in  
[1,+\infty[$,  
\begin{eqnarray*}A(x,x')&\leq& C_p\,\sup_{t\in [0,T]}\, E \Big|\int_0^t ds  
\int_Q dy\int_Q |H(t,x;s,y)-H(t,x';s,y)|\, |u(s,y)|\\  
&&\qquad\qquad \qquad  \times\, f(y-z)\, |H(t,x;s,z)-H(t,x';s,z)|\,|u(s,z)|\,  
 dz\Big|^p\, .  
\end{eqnarray*}  
We prove that for $\Delta(x,x')= \sup_{t\in [0;T]}
\int_0^t ds \int_Q dy\int_Q dz \, |H(t,x;s,y)-H(t,x';s,y)|\,  
f(y-z)$\linebreak[3] $\times |H(t,x;s,z)-H(t,x';s,z)|$, one has 
\begin{equation}\label{Deltax}\Delta(x,x')\leq C_p\, |x-x'|^a\, ,  
\end{equation}  
 provided that $H$ satisfies either (\ref{covISxd=}) or (\ref{covISxdif}).  
  Assuming that  (\ref{Deltax}) holds, H\"older's and Schwarz's inequalities  
   yield  
\begin{eqnarray*}  
A(x,x')&\leq&C_p\; \Delta(x,x')^{p-1} \, \sup_{t\in [0,T]} \,  
\int_0^tds \int_Q dy\,  |H(t,x;s,y)-H(t,x';s,y)|\, \\  
&&\quad\times\int_Q f(y-z)\, |H(t,x;s,z)-H(t,x';s,z)|\,  
E\left(|u(s,y)|^{p}\, |u(s,z)|^p\right) dz\\ &\leq&C_p\;  
\Delta(x,x')^{p} \,  \sup_{(s,y)\in [0,T]\times Q}  
E(|u(s,y)|^{2p})\leq C_p\; |x-x'|^{ap}\, .  
\end{eqnarray*}  
In order to prove (\ref{Deltax}), we use  Taylor's formula, the convexity 
of $Q$  and the  
inequalities (\ref{majoH}) and (\ref{majoDH}); thus for any $a\in  
]0,1[$,  $\Delta(x,x')\leq |x-x'|^a\, (T_1+T_2)\, $, where  
\begin{eqnarray*}  
T_1&=&\sup_{x\in Q}\int_0^T s^{-(2\alpha +a\delta)}\, ds\int_Q dy \int_Q  
\exp\left(-c\frac{|x-y|^\beta}{s^\gamma}\right)\,  f(y-z)\,  
\exp\left(-c\frac{|x-z|^\beta}{s^\gamma}\right)\, dz\, ,\\  
T_2&=&\int_0^T s^{-(2\alpha +a\delta)}\, ds\int_Q dy \int_Q  
\exp\left(-c\frac{|x-y|^\beta}{s^\gamma}\right)\,  f(y-z)\,  
\exp\left(-c\frac{|x'-z|^\beta}{s^\gamma}\right)\, dz\, .  
\end{eqnarray*}  
Replacing $2\alpha$ by $2\alpha +a\delta$, the arguments used to  
prove part (i) of Lemma \ref{CNS} show that if either  
(\ref{covISxd=}) or (\ref{covISxdif}) holds, then $T_1<+\infty$.  
To study $T_2$, we have to distinguish several cases. Let  
$0<c_2<1$,  $k\geq 1$ be such that $c_2^k <\frac{1}{3}$ and set  
$\varepsilon = c_2^{k+1}$,  
$\overline{\varepsilon}=\frac{c_2^k}{1-c_2^{k+1}}$; then  
$\frac{1}{1+\varepsilon}>\overline{\varepsilon}$ and $\varepsilon  
<1$ . We study three cases:  
\smallskip

{\em Case 1} If $|x-y|\geq \varepsilon \, |x'-y|$ or $|x'-z|\geq  
\varepsilon \, |x-z|$, we have (changing the constant in the  
exponential functions) $T_2\leq C(\varepsilon)\, T_1$, and the  
proof is complete.  
  
{\em Case 2} If $|x-y|<\varepsilon \, |x'-y|$,  $|x'-z|<  
\varepsilon \, |x-z|$ and $|y-z|\leq \overline{\varepsilon}\,  
|x-x'|$, then $|x-z|\geq \frac{|x-x'|}{1+\varepsilon}$ and we have  
\[|x-y|\geq \big|\, |x-z|-|y-z|\, \big| \geq \left(\frac{1}{1+\varepsilon}-  
\overline{\varepsilon}\right)\,  |x-x'| \, \geq \left(\frac{1}{  
\bar{\varepsilon}(1+\varepsilon)}- 1\right)\, |y-z|\, .\]  
This implies for $\bar{c}=c \min\Big(1,\frac{1}{\overline{\varepsilon}  
(1+\varepsilon )} -1 \Big)$  which is positive by the choice of $k$:  
\[  
T_2\leq\int_0^T s^{-(2\alpha +a\delta)} \, ds\int_Q dy \int_Q \,  
\exp\left(-\bar{c}\frac{|y-z|^\beta}{s^\gamma}\right)\,  f(y-z)\,  
\exp\left(-\bar{c}\frac{|x'-z|^\beta}{s^\gamma}\right)\, dz\, ,  
\]  
and since this is the upper estimate of (\ref{numero}), the proof is again  
 concluded by an argument similar to that in Lemma \ref{CNS} (i), with  
 $2\alpha+a\delta$ instead of $2\alpha$.  
 \smallskip

{\em Case 3} Suppose finally that $|x-y|< \varepsilon \, |x'-y|$,  
$|x'-z|< \varepsilon \, |x-z|$ and $|y-z|>  
\overline{\varepsilon}|x-x'|$. Then $|x'-y|\leq  
\frac{|x-x'|}{1-\varepsilon}$ and $ |x-y|\leq  
\frac{\varepsilon}{1-\varepsilon}\, |x-x'|\leq  
 \frac{\varepsilon}{(1-\varepsilon)\, \overline{\varepsilon}}\,  
 |y-z|=c_2\, |y-z|\, ,$  
so that $\frac{c_1}{c_2}\, |x-y|\leq  c_1\, |y-z|$. Since $f$  
satisfies (\ref{fmonoton}), we have $f(y-z)\leq C_1\,  
f\left(\frac{c_1}{c_2}\, (y-x)\right)$. Set  
$\bar{y}=x+\frac{c_1}{c_2}\, (y-x)$ and $\tilde{c}=c\, \min  
\Big(\frac{c_1}{c_2}, 1\Big)$; then  
%\begin{eqnarray*}  
\[  
T_2\leq C\, \int_0^T s^{-(2\alpha + a\delta)} \, ds% \\  
%&&\qquad\qquad\qquad \times  
\int_Q d\bar{y} \int_Q \,  
\exp\left(-\tilde{c}\frac{|x-\bar{y}|^\beta}{s^\gamma}\right)\,  
f(\bar{y}-x)\,  
\exp\left(-\tilde{c}\frac{|x'-z|^\beta}{s^\gamma}\right)\, dz\, ,  
\]  
%\end{eqnarray*}  
and again the proof is complete, since the right hand-side is  
similar to (\ref{numero}).\smallskip

(ii) For $0\leq t<t'\leq T$, set  
\begin{eqnarray*}  
\overline{\Delta}(t,t')&=& \sup_{x\in Q} \Bigg|\int_0^t ds \int_Q  
\, dy \; |H(t,x;s,y)-H(t',x;s,y)|\\ &&\qquad\qquad\times \int_Q  
f(y-z)\, |H(t,x;s,z)-H(t',x;s,z)|\, dz\, \Bigg|^p\, .  
\end{eqnarray*}  
Again we prove that under either condition (\ref{covIStd=}) or  
(\ref{covIStdif}) we have  
\begin{eqnarray}\label{Deltat}  
\overline{\Delta}(t,t')\leq C |t-t'|^b\, .  
\end{eqnarray}  
If (\ref{Deltat}) holds, using again  
Burkholder's, H\"older's and Schwarz's inequalities, we deduce that  
%\begin{eqnarray*}  
\[B(t,t')\leq C_p\, \overline{\Delta}(t,t')^p \, \sup_{(t,y)\in  
[0,T]\times Q} E(|u(s,y)|^{2p})%\\ &\leq&  
\leq C_p\, |t-t'|^{bp}\, .\]  
%\end{eqnarray*}  
We now prove (\ref{Deltat}). Using (\ref{majoH}), (\ref{majoDH})  
and Taylor's formula, we obtain for $h=t'-t$: for any $b\in]0,1[$,  
$ \overline{\Delta}(t,t')\leq  |t-t'|^{b}\, (T'_1 + T'_2 )$,  where  
\begin{eqnarray*}  
T'_1&=&\sup_{x\in Q} \int_0^t s^{-(2\alpha +b\eta )}\, ds \int_Q  
dy\int_Q \exp\left(-c\, \frac{|x-y|^\beta}{s^\gamma}\right)\,  
f(y-z)\, \exp \left(-c\, \frac{|x-z|^\beta}{s^\gamma}\right)\, dz  
\, ,\\ 
T'_2&=&\sup_{x\in Q} \int_0^t s^{-(\alpha +b\eta )}\,  
(s+h)^{-\alpha}\,  ds \int_Q dy  \exp\left( -c\,  
\frac{|x-y|^{\beta} }{s^{\gamma}}\right) \\  
 & & \qquad \qquad \times \, \int_Q f(y-z)\, \exp\left( -c\,  
 \frac{|x-z|^\beta}{(s+h)^\gamma}\right)\, dz \, .  
\end{eqnarray*}  
Clearly, $T'_1$ is similar to $T_1$ with $b\eta $   instead of  
$a\delta$; thus the proof of (i) yields  
(\ref{acroiISt}) if either (\ref{covIStd=}) or (\ref{covIStdif})  
holds. To estimate $T'_2$, we distinguish two cases.\smallskip

{\em Case 1} If $|x-y|\leq c_1\, |y-z|$, condition (\ref{fmonoton}) on $f$
 yields 
\begin{eqnarray*}  
 T'_2&\leq& C_1\, \int_0^t s^{-(\alpha +b\eta )}\, (s+h)^{-\alpha}\, ds  
 \int_Q \exp\left( -c\, \frac{|x-y|^{\beta}}{s^{\gamma}}\right)\, f(x-y)  
 \, dy \\  
 & &\qquad\qquad \times \int_Q \,  
 \exp\left( -c\, \frac{|x-z|^\beta}{(s+h)^{\gamma}}\right)\, dz \, .  
 \end{eqnarray*}  
Set  $v=x-y$; then since $\alpha \geq \frac{\gamma}{\beta}d$, we  
have for some $R>0$,  
\begin{eqnarray*}  
T'_2&\leq & C_1\, \int_0^t s^{-(\alpha +b\eta )}\,  
(s+h)^{-\alpha+\frac{\gamma}{\beta} d}\, ds \int_{B_d(0,R)}  
\exp\left(-c\, \frac{|v|^\beta}{s^\gamma}\right)\, f(v)\, dv\\  
&\leq&C_1\, \int_0^t s^{-(2\alpha +b\eta )+\frac{\gamma}{\beta}  
d}\, ds \int_{B_d(0,R)} \exp\left(-c\,  
\frac{|v|^\beta}{s^\gamma}\right)\, f(v)\, dv\, .  
\end{eqnarray*}  
This last upper estimate is similar to the right hand side of (\ref{inter})  
 with $2\alpha +b$ instead of $2\alpha$; thus the end of the proof of  
  Lemma \ref{CNS} (i) concludes the proof.\smallskip

{\em Case 2} If $|x-y|>c_1\, |y-z|$, then for $\bar{c}=\min(c, c\,  
 c_1^\beta)$, we have  
\begin{eqnarray*} T'_2&\leq &\sup_{x\in Q}\int_0^t s^{-(\alpha +b\eta )}\,  
(s+h)^{-\alpha}\, ds \int_Q dy \\  
&&\qquad\qquad \times \int_Q  
\exp\left(-\bar{c}\, \frac{|y-z|^\beta}{s^\gamma}\right)\, f(y-z)\,  
 \exp\left(-\bar{c}\, \frac{|x-z|^\beta}{(s+h)^\gamma}\right)\, dz\, .  
\end{eqnarray*}  
Since $(s+h)^{-\alpha+\frac{\gamma}{\beta}d} \leq  
s^{-\alpha+\frac{\gamma}{\beta}d}$, the change of variables $\zeta  
= x-z$ and $v=y-z$ shows that $T'_2$ is dominated by the right  
hand side of (\ref{inter}) with $2\alpha +b$ instead of  
$2\alpha$ and $\bar{c}$ instead of $c$; this concludes the proof  
of (\ref{acroiISt})  
\smallskip

Finally, using again Burkholder's inequality and (\ref{majoH}), we  
have ${\displaystyle C(t,t')\leq C_p \, T_3^p\, \sup_{(s,y)\in  
[0,T]\times Q}}$\linebreak[3] $E(|u(s,y)|^{2p})$,  where  
\[T_3=\sup_{x\in Q} \int_0^{t'-t} s^{-2\alpha}\, ds \int_Q  dy\int_Q \,  
\exp\left(-c\, \frac{|x-y|^\beta}{s^\gamma}\right)\, f(y-z)\,  
\exp\left(-c\, \frac{|x-z|^\beta}{s^\gamma}\right)\, dz\, .  
\]  
Computations similar to those in the proof of Lemma \ref{CNS} (i)  
imply that for some $R>0$, $ T_3\leq \int_0^{t'-t}  
%s^{-2\alpha+\frac{\gamma}{\beta} d}\, ds \int_{B_d(0,R)}  
%\exp\left(-\frac{|v|^\beta}{s^\gamma}\right)\, f(v)\, dv\,  
\psi(s)\, ds$, where $\psi$ is defined by (\ref{defpsi}). Fubini's  
theorem and H\"older's  
 inequality with respect to $ds$ with the conjugate exponents $\lambda =  
 (b\, \eta)^{-1}$ and $\mu$ imply  
\[T_3\leq |t'-t|^{b} \int_{B_d(0,R)}  f(v)\, dv \left( \int_0^T  
\exp\left(-c\, \frac{|v|^\beta}{s^\gamma}\right)\,  
s^{\mu(\frac{\gamma}{\beta}d-2\alpha)}\,  
ds\right)^{\frac{1}{\mu}}\, .\] For $v\neq 0$ set $r=|v|^{\beta}\,  
s^{-\gamma}$; then since $\frac{1}{\mu}=1-b$, we obtain  
\[T_3\leq C\, |t'-t|^{b}\int_{B_d(0,R)} f(v) \,  
|v|^{d+\frac{\beta}{\gamma}(1-2\alpha-b)}\, dv \left(  
\int_{|v|^\beta\, T^{-\gamma} }^{+\infty} r^{-1+\frac{\mu}{\gamma}  
[2\alpha -\frac{\gamma}{\beta}d -1+b ]}\, \exp(-c\, r)\,  
dr\right)^{\frac{1}\mu} \, .\]  
 As in the proof of Lemma \ref{CNS}, we distinguish three cases,  
according to the power of $r$ in the last integral, with $2\alpha  
+b$ instead of $2\alpha$; this concludes the proof of  
(\ref{IStt'}).  
\epr

\section{Existence of solutions}  
\subsection{The case of Lipschitz coefficients}  
  
Let $\sigma : [0,T]\times Q\times \mathbb{R}\rightarrow  
\mathbb{R}$ and $b_i : [0,T]\times Q\times \mathbb{R}\rightarrow  
\mathbb{R}$, $1\leq i\leq N$ be continuous functions such that the  
following boundness and Lipschitz conditions hold:  
\begin{description}  
\item %{\makebox[24pt]  
\textbf{(L1) }%}  
  Uniform linear growth with respect to  
the last variable:  
 for every $y\in Q$  
\begin{equation}\label{bound}\sup_{(t,x)\in [0,T]\times Q}  
\Big( |\sigma(t,x,y)|+\sum_{i=1}^N|b_i(t,x,y)|\Big) \leq  
C(1+|y|)\, .  
\end{equation}  
\item %{\makebox[24pt]  
\textbf{(L2) }%}  
 Uniform Lipschitz condition with respect  
 to the last variable: for any $y,z\in Q$  
\begin{equation}\label{lip}  
\sup_{(t,x)\in [0,T]\times Q} \Big(|\sigma(t,x,y)-\sigma(t,x,z)|+  
\sum_{i=1}^N |b_i(t,x,y)-b_i(t,x,z)|\Big)\leq C\, |y-z|\, .  
\end{equation}  
\end{description}  
We then consider the following non-linear evolution equation for  
$t\in [0,T]$  
 and $x\in Q$:  
\begin{eqnarray}\nonumber  
u(t,x)&=&\int_Q G(t,x;0,y)\,  u_0(y)\, dy + \int_0^t\int_Q  
G(t,x;s,y)\, \sigma(s,y,u(s,y))\, F(ds,dy)\\ &&\qquad +  
\sum_{i=1}^N \int_0^t\int_Q H_i(t,x;s,y)\, b_i(s,y,u(s,y))\, dy \,  
ds \, , \label{evolip}  
\end{eqnarray}  
 for a  function $u_0 : Q\rightarrow \mathbb{R}\in L^2(Q)$.  
\smallskip

In this section we will make the following assumptions  
(restricting ourselves to the case $G$ is the Green function of  
an operator which is parabolic in the sense of Petrovski\u{\i},  
 and the functions $H_i$ are  partial derivatives of $G$ with respect to the  
  space variable $y$):  
\begin{description}  
\item %{\makebox[24pt]  
\textbf{(C2) }%}  
 The continuous function $G : ([0,T]\times Q)^2\rightarrow\mathbb{R}$  
(respectively for $i=1,\cdots ,N$ each  continuous functions $H_i  
:([0,T]\times Q)^2\rightarrow\mathbb{R}$)   satisfies (\ref{majoH})  
with constants $\beta, \gamma$ and $\alpha =\frac{\gamma}{\beta}d$  
(respectively $\alpha_i,\beta$ and $\gamma$).  
  
\noindent\item% {\makebox[24pt]  
\textbf{(C3) }%}  
 The covariance of the  
Gaussian process  
 $F$ is defined in terms of $f$ by (\ref{cov}), and the constants $\alpha,  
 \beta,\gamma$ in (C2) satisfy either  
\begin{equation}\label{covd=}  
\int_{B_d(0,1)} f(v)\, \, \ln\big(|v|^{-1}\big)\, dv <\infty,  
\quad if \quad d=\frac{\beta}{\gamma}\, ,  
\end{equation}  
or  
\begin{equation}\label{covdiff}  
\int_{B_d(0,1)} f(v)\, |v|^{-[d-\frac{\beta}{\gamma}]^+}\,  dv <\infty,  
\quad if \quad  d\neq \frac{\beta}{\gamma}\, .  
\end{equation}  
\item %{\makebox[24pt]  
\textbf{(C4) }%}  
 The constants $\alpha_i,\beta,\gamma$  
in (C2) satisfy $\alpha_i<\alpha +1$ for every $i\in \{1, \cdots,  
N\}$.  
\end{description}  
Condition  (C3) will allow to define stochastic integrals of $G(t,x;.)$ with  
 respect to the noise $F$, while  (C4) will allow to define deterministic  
integrals involving $H_i(t,x;.)$.  
  
We at first study moment estimates of deterministic integrals. Let  
$\lambda, \rho\in [1,+\infty]$; for $v\in L^{\lambda}([0,T],  
L^{\rho}(Q))$, and $0\leq t_0\leq t\leq T$, $x\in Q$ set  
\begin{equation}\label{J}  
J(v)(t_0,t,x)=\int_{t_0}^t\int_Q H(t,x;s,y)\, v(s,y) dy ds\,.  
\end{equation}  
The following lemma provides $L^q$ estimates of $J(v)(t_0,t,.)$ in terms of  
$L^{\rho}$ estimates of $v(s,.)$. It extends similar results proved in  
I.~Gy\"ongy \cite{Gy} and C.~Cardon-Weber \cite{CW1}.  
\begin{lem}.\label{convolution}  
Fix $\rho\in [1,+\infty]$, $q\in [\rho,+\infty]$, and let $r$ be  
defined by $\frac{1}{r}=\frac{1}{q}-\frac{1}{\rho}+1$. Then for  
$0\leq t_0\leq t\leq T$  
\begin{equation}\label{momentqJ(v)}  
\|J(v)(t_0,t,.)\|_q\leq C\int_{t_0}^t (t-s)^{-\alpha +  
\frac{\gamma d}{\beta r}}\, \|v(s,.)\|_{\rho}\, ds\, .  
\end{equation}  
Hence given any $\lambda\in [0,+\infty]$, if $ \alpha  
+\frac{1}{\lambda}< \frac{\gamma d}{\beta r}+1$, then  $J(0,.,.)$  
is a bounded operator from $L^{\lambda}([0,T], L^{\rho}(Q))$ into  
$L^{\infty}([0,T], L^q(Q))$.  
\end{lem}  
\bpr Using Minkowski's inequality, (\ref{majoH}),  then Young's  
inequality with $\frac{1}{q}=\frac{1}{r}+\frac{1}{\rho}-1$ and  
(\ref{intDH}), we obtain  
\begin{eqnarray*}  
\| J(v)(t_0,t,.)\|_q&\leq&  
 C\int_{t_0}^t\left\| \int_Q (t-s)^{-\alpha}\, \exp\left( -c  
 \frac{|.-y|^\beta}{(t-s)^{\gamma}}\right)\, |v(s,y)|\, dy\right\|_q\, ds\\  
&\leq & C\int_{t_0}^t (t-s)^{-\alpha}\, \|v(s,.)\|_{\rho} \, \left\|  
\exp\left( -c\frac{|.|^\beta}{(t-s)^{\gamma}}\right)\right\|_r\, ds\\  
&\leq &C\int_{t_0}^t(t-s)^{-\alpha +\frac{\gamma d}{\beta r}}\,  
\|v(s,.)\|_{\rho}\, ds\, .  
\end{eqnarray*}  
Finally, H\"{o}lder's inequality applied with $\lambda$ and $\mu=  
\frac{\lambda}{\lambda -1}$ yields that  
\[ \| J(v)(t_0,t,.)\|_q\leq C\, \left(\int_{t_0}^t (t-s)^{\mu (-\alpha +  
\frac{\gamma d}{\beta r})}\, ds\right)^{\frac{1}{\mu}} \left(\int_{t_0}^t  
 \|v(s,.)\|_{\rho}^{\lambda}\, ds\right)^{\frac{1}{\lambda}}\, .\]  
Hence for $t_0=0$, the right hand side is finite and bounded with  
respect to $t\in [0,T]$  if and only if $\mu (-\alpha  
+\frac{\gamma d}{\beta r}) >-1$; this completes the proof.  
\hfill $\Box$%\epr  
\medskip

The following result proves that the evolution equation (\ref{evolip}) has a  
 unique solution with moments of all finite order.  
However, in order to prove that, when $u_0\in L^q(Q)$ for $2\leq  
q<+\infty$, the $\|\;\|_q$-norm of the solution has bounded $L^p$  
moments for $q<p<+\infty$, we have to reinforce condition (C3) as  
follows (clearly when $p=q$, the conditions (C3) and (C'3)(q,p)  
coincide, while if $p<q$, (C'3)(q,p) implies (C3)):  
\medskip

 {\rm  
\textbf{(C'3)(q,p)}  Let $f$ define the covariance of the Gaussian  
noise according to (\ref{cov}),  $2\leq q\leq p<+\infty$; the  
constants $\alpha, \beta,\gamma$ in (C2) satisfy one of the  
following conditions:  
\begin{equation}\label{covd=b}  
\int_{B_d(0,1)} f(v)\, \, \ln\big(|v|^{-1}\big)\, dv <\infty,  
\quad if \quad \frac{\beta}{\gamma}(2\alpha-1)=\frac{q}{p}d\, ,  
\end{equation}  
or  
\begin{equation}\label{covdiffb}  
\int_{B_d(0,1)} f(v)\,  
|v|^{-\left[\frac{\beta}{\gamma}(2\alpha-1)-d\frac{q}{p}\right]^+}\, dv  
<\infty, \quad if \quad \frac{\beta}{\gamma}(2\alpha-1)\neq\frac{q}{p}d\, .  
\end{equation}  
\begin{Th}.\label{existlip}  
Suppose that the functions $G$ and $H_i$ , $1\leq i\leq N$ satisfy  
the conditions  (C2) and (C4) and that the functions $\sigma$ and  
$b_i$, $1\leq i\leq N$ satisfy the assumptions  (L1) and  (L2).

(i) Let $u_0\in L^\infty(Q)$, and let $F$ denote either the space-time white noise
if $\alpha <1$, or a Gaussian process with covariance defined by  
(\ref{cov})  such that (C3) holds. Then the  
 evolution equation (\ref{evolip}) has a unique solution $u\in L^{\infty}  
 ([0,T], L^\infty(Q))$, such that for any $p\in [1,\infty[$, 
 \begin{equation}\label{momlip}  
\sup_{(t,x)\in [0,T]\times Q} E\big( |u(t,x)|^p\big) < +\infty\, .
%\quad if\quad   q=+\infty \quad and\quad 1\leq p<+\infty ,  
\end{equation}  

(ii) Let $u_0\in L^q(Q)$ for $2\leq q<+\infty$, let $p\in [q,+\infty[$.
Suppose that the following assumptions holds:

\qquad (a) $F$ is the space-time white noise, $\alpha <1$ and $p<\frac{2\alpha}{\left[ \frac{2\alpha}{ q}-1+\alpha\right]^+}$.

\qquad (b) $F$ is a Gaussian process with covariance defined by  
(\ref{cov})  such that (C'3)(q,p) holds.

\noindent Then the  evolution equation (\ref{evolip}) has a unique solution $u\in L^{\infty}  
 ([0,T], L^q(Q))$, such that 
 \begin{equation}\label{momlipq}  
\sup_{t\in [0,T]} E\big( \|u(t,.)\|_q^p\big) < +\infty .  
\end{equation}  
\end{Th}

% and  that  $F$ is a Gaussian process with covariance defined by  
%(\ref{cov}),  such that (C3) holds.\smallskip
%Then for any  initial condition $u_0\in L^q(Q)$, $2\leq q\leq +\infty$,  the  
 %evolution equation (\ref{evolip}) has a unique solution $u\in L^{\infty}  
% ([0,T], L^q(Q))$. Furthermore,  
%while if  $2\leq q\leq p <+\infty$, and  if $f$ satisfies  
%condition (C'3)(q,p), then  
%\begin{equation}\label{momlipq}  
%\sup_{t\in [0,T]} E\big( \|u(t,.)\|_q^p\big) < +\infty .  
%\end{equation}  
%\end{Th}  

\bpr  In the case of the space-time white noise, the proof which is easier
and more classical is omitted, except that of (\ref{momlipq}) in case (ii).
Unless specified otherwise, we assume that $F$ is Gaussian with a space-correlation 
function $f$. We use the following Picard iteration scheme;
 $ u_0(t,x)=  
G_t u_0(x)=\int_Q G(t,x;0,y)\,  u_0(y)\, dy $ and  for $n> 0$ let  
\begin{eqnarray}  
\nonumber u_{n+1}(t,x)&=&u_0(t,x) + \int_0^t\int_Q G(t,x;s,y)\,  
\sigma(s,y,u_n(s,y))\, F(ds,dy)\\ \label{picardlip} &&\qquad +  
\sum_{i=1}^N \int_0^t\int_Q H_i(t,x;s,y)\, b_i(s,y,u_n(s,y))\, dy  
\, ds \, . 
\end{eqnarray}

{\em Case (ii)} Let $2\leq q\leq p<+\infty$ and suppose that condition (b) 
holds;   set  
\[M_n(t)=E\Big(\|u_n(t,.)\|_q^p \Big)\, ,\]  
and let $\psi_p$ be  the function defined by ${\displaystyle  
\psi_p(t)=t^{\alpha \left(-2+\frac{q}{p}\right)}  
\int_{B_d(0,R)}\exp\left( -c\frac{\vert  
v\vert^\beta}{t^\gamma}\right)f(v) dv\, .}$  Using (C'3)(q,p),  
computations similar to that proving (\ref{intconv}) from  
(\ref{inter}) using (\ref{defpsi}) show that  
 $\psi_p$ is integrable; set $I_p=\int_0^T \psi_p(s) ds < \infty$  
and let  
\[\varphi_p(s)=C_p\, \Big( \psi_p(s) + \sum_{i=1}^N s^{-\alpha_i +  
\alpha} \Big)\, ;\] the assumptions (C'3)(q,p), (C4),  
(\ref{intDH}) and the proof of lemma \ref{CNS} imply that  
$\varphi_p \in L^1_+ ([0,T])$. Let $\psi$ and $I$ be defined as in  
the proof of lemma \ref{CNS}; then for $q\leq p$, $\psi_p\geq  
\psi_q=\psi$ and $I_q\geq I_p=I$. We prove that  
\begin{equation}\label{0}  
\sup_{t\in[0,T]} M_0(t)<+\infty\, ,  
\end{equation}  
and for any $n\geq 0$, $ t\in[0,T]$,  
\begin{equation}\label{dalang}  
M_{n+1}(t)\leq \int_0^t\varphi_p(t-s)\, [1+M_n(s)]\, ds\, ;  
\end{equation}  
then Lemma 15 in \cite{D} shows that  
\begin{equation}\label{majq}  
\sup_n\, \sup_{0\leq t\leq T} E\left(\|u_n(t,.)\|_q^p\right) <+\infty\, .  
\end{equation}  
Fubini's theorem, (\ref{majoH}), H\"older's inequality and (\ref{intDH})  
imply that  
\begin{eqnarray*}  
\sup_{0\leq t\leq T}M_0(t)&\leq& C\,\sup_{0\leq t\leq T}\left [  
\int_Q\left|\int_Q u_0(y)\, t^{-\alpha }\, \exp\left(-c\,  
\frac{|x-y|^\beta}{t^\gamma}\right)\, dy\, \right|^q\, dx \right  
]^{\frac{p}{q}}\\ &\leq&C \left [ \int_Q |u_0(y)|^q\, \left(  
t^{-\alpha }\, \int_Q\exp\left(-\bar{c}\,  \frac{|x-y|^\beta}{  
t^\gamma}\right)\, dx\, \right)\, dy\,\right ]^{\frac{p}{q}} \leq  
C\, \|u_0\|_q^p\,  .  
\end{eqnarray*}  
  
We now prove (\ref{dalang}); for $n\geq 0$, $M_{n+1}(t)\leq  
C_p\left[M_0(t)+T^1_n(t,p)+\sum_{i=1}^NT^2_{n,i}(t,p)\right]$,  
where for $q\leq p<+\infty$,  
\begin{eqnarray*}  
T_n^1(t,p)&=&E\left(\left\|\int_0^t\int_Q G(t,.;s,y)\, \sigma(s,y,u_n(s,y))  
\, F(ds,dy)\right\|_q^p \right)\, ,\\  
T_n^{2,i}(t,p)&=&E\left(\left\|\int_0^t\int_Q H_i(t,.;s,y)\,  
 b_i(s,y,u_n(s,y))\, dy\, ds\right\|_q^p \right)  \, .
\end{eqnarray*}  
Since $\beta\geq 1$, $|y-z|^\beta\leq 2^{\beta-1}\, \big[  
 |x-y|^\beta+|x-z|^\beta\big]$. Therefore,  Fubini's theorem,  Burkholder's  
 and H\"older's inequalities yield the existence of a constant $\bar{c}$  
 such  that  
\begin{eqnarray*}  
T_n^1(t,p)&\leq&C_p\, \int_Q dx \, E\left( \left|  \int_0^t ds  
(t-s)^{-2\alpha} \int_Q dy \,  
 \exp\left( -\bar{c}\frac{|x-y|^\beta}{(t-s)^\gamma}\right)\,  
\big(1+|u_n(s,y)|\big)\right.\right.\\ &&\qquad\qquad\qquad  
\left.\left. \times\int_Q f(y-z)\,  
\exp\left(-\bar{c}\frac{|x-y|^\beta}{(t-s)^\gamma}\right)\,  
\big(1+|u_n(s,z)|\big)\, dz \right|^{\frac{p}{2}} \right)\\  
&\leq&C_p\, I_q^{\frac{p}{2}-\frac{p}{q}}\, \int_Q dx\, E\left(  
\left[ \int_0^t (t-s)^{-2\alpha }\int_Q dy\, \exp\left( -\bar{c}\,  
\frac{|x-y|^\beta}{(t-s)^\gamma}\right)\,  
\left(1+|u_n(s,y)|^{\frac{q}{2}}\right)\right.\right .\\  
&&\qquad\qquad\qquad\times\left.\left .\int_Q f(y-z)\, \exp\left(  
-\bar{c}\, \frac{|y-z|^\beta}{(t-s)^\gamma}\right)\,  
\left(1+|u_n(s,z)|^{\frac{q}{2}}\right)\, dz  
ds\right]^{\frac{p}{q}} \right)\, .  
\end{eqnarray*}  
Fubini's theorem and Jensen's inequality imply that  
  \begin{eqnarray*}  
T_n^1(t,p)&\leq&C_p\, I_q^{\frac{p}{2}-\frac{p}{q}}\, E\Bigg(  
\Bigg[ \int_0^t (t-s)^{-2\alpha }\Big \Vert \int_Q dy\, \exp\left(  
-\bar{c}\, \frac{|.-y|^\beta}{(t-s)^\gamma}\right)\,  
\left(1+|u_n(s,y)|^{\frac{q}{2}}\right)\\  
&&\qquad\qquad\qquad\times\int_Q f(y-z)\, \exp\left( -\bar{c}\,  
\frac{|y-z|^\beta}{(t-s)^\gamma}\right)\,  
\left(1+|u_n(s,z)|^{\frac{q}{2}}\right)\,  
dz\Big\Vert_{\frac{p}{q}}ds \Bigg]^{\frac{p}{q}} \Bigg)\, .  
\end{eqnarray*}  
Then Young's, Schwarz's   inequalities, (\ref{intDH}), H\"older's  
inequality with respect to $\psi_p(t-s)\, ds$ and $I_q\leq I_p$  
yield  
\begin{eqnarray}  
T_n^1(t,p)&\leq&C_p\, I_q^{\frac{p}{2}-\frac{p}{q}}\,  
E\Bigg(\Bigg\vert \int_0^t (t-s)^{-2\alpha}\, \Big\Vert  
 \exp\Big(-\bar{c}\frac{|.|^\beta}{(t-s)^\gamma}\Big)\, *\,  
\Big[ \, \left( 1+|u_n(s,.)|^{\frac{q}{2}} \right) \nonumber\\  
&&\qquad\qquad\times\Big\{ \Big( f(.)\,  
\exp\big(-\bar{c}\frac{|.|^\beta}{(t-s)^\gamma}\big)\Big)\, *\,  
\Big(1+|u_n(s,.)|^{\frac{q}{2}}\Big)\Big\} \Big]\,  
\Big\Vert_{L^{\frac{p}{q}}(Q,dx)}\,ds \Bigg\vert^{\frac{p}{q}}  
\Bigg)\nonumber \\  
 &\leq&C_p\, I_q^{\frac{p}{2}-\frac{p}{q}} \,  
 E\Bigg(\Bigg\vert \int_0^t (t-s)^{-2\alpha}\,  
\Big\|\exp\Big(-\bar{c}\, \frac{|.|^\beta}{(t-s)^\gamma }\Big)  
\Big\|_{L^{\frac{p}{q}}(Q,dx)} \,\left(  
1+\|u_n(s,.)\|_q^{\frac{q}{2}}\right)\nonumber\\ &&\qquad\qquad  
\times \Big\|f(.)\, \exp\Big(-\bar{c}  
\frac{|.|^\beta}{(t-s)^\gamma}\Big)\Big\|_{L^1(Q,dx)}\,  
\left(1+\|u_n(s,.)\|_q^{\frac{q}{2}}\right)\,  
ds\Bigg\vert^{\frac{p}{q}} \Bigg) \nonumber\\  
 &\leq& C_p\,  
I_q^{\frac{p}{2}-\frac{p}{q}} \,E\left (\left \vert \int_0^t  
 (t-s)^{\alpha \left(-2+\frac{q}{p}\right)}\,  \left[1+\left(  
 \|u_n(s,.)\|_q^q\right)\,\right]  
\right.\right.\nonumber\\ & & \qquad \qquad \left. \left. \times\,  
\int_Q f(v)\,  
\exp\left(-\bar{c}\frac{|v|^\beta}{(t-s)^\gamma}\right)\, dv  
ds\right \vert^{\frac{p}{q}}\right)\nonumber \\ &\leq& C_q\,  
I_p^{\frac{p}{2}-1}\int_0^t\psi_p(t-s)\, \left[1+M_n(s)\right]\,  
ds\, .\label{majoTn1}  
\end{eqnarray}  
For every $1\leq i\leq N$, using (C4), (L1), (\ref{momentqJ(v)})  
with $\rho=q$, $r=1$, and H\"older's inequality (since  
$\alpha_i<\alpha +1$) and Fubini's theorem, we deduce that for  
$q\leq p<+\infty$,  
\begin{eqnarray*}  
 T_{n,i}^2(t,p)&\leq&C\, E\left(\left|\int_0^t(t-s)^{-\alpha_i +  
\alpha}\, \left(1+\|u_n(s,.)\|_q\right)\, ds \right|^p \right)\\  
&\leq&C\, \int_0^t(t-s)^{-\alpha_i +\alpha}\,  
 \left[1+E( \|u_n(s,.)\|_q^p)\, \right]\, ds\, .  
\end{eqnarray*}  
This concludes the proof of (\ref{dalang}). Let  
$\Delta_n(t)=E\left(\left\|u_{n+1}(t,.)-u_n(t,.)\right\|_q^p\right)$;  
a similar computation using the global Lipschitz property (L2) of  
the coefficients with respect to the last variable shows that  
\begin{equation}\label{delta}  
\Delta_{n+1}(t)\leq C_p\,  \int_0^t \varphi_p(t-s)\, \Delta_n(s)\, ds  
\end{equation}  
where the function $\varphi_p$ is the previous one.  Using again  
Lemma 15 in \cite{D}, we conclude that $\sum_{n\geq 0}  
\Delta_n(t)$ converges uniformly on $[0,T]$. Therefore, usual  
arguments show that the solution $u$ to (\ref{evolip}) exists in  
$L^{\infty}([0,T],L^q(Q))$ and satisfies (\ref{momlipq}).  
\bigskip

We now suppose that condition (a) holds. 
Set $M_n(t)=E\left(\|u_n(t,.)\|_q^p\right)$.
 According to the results proved above, it suffices to check that (using the
  previous notations), $T_n^1(t,p)\leq C_p\, \int_0^t (t-s)^{-a}\, (1+
  M_n(s,p))\, ds$ for some $a<1$.
 Using H\"older's and Burkholder's inequalities, Fubini's theorem, then
 (\ref{momentqJ(v)}) with $\frac{p}{2}$
 instead of $q$, $2\alpha$ instead of $\alpha$ and $\rho=\frac{q}{2}$ , we
 deduce that for $1+\frac{2}{p}=\frac{2}{q}+\frac{1}{r}$, and $a=(2-
 \frac{1}{r})\alpha <1$  by the choice of $p$, we have
 \begin{eqnarray*}
 T_n^1(t,p)&\leq& C_p\, E\left(\int_Q  \left| \int_0^t \int_Q
 (t-s)^{-2\alpha}\,
 \exp\left(-c\frac{|x-y|^\beta}{(t-s)^\gamma}\right)\, (1+|u_n(s,y)|^2)\,
 dy\, ds\right|^{\frac{p}{2}}\, dx\right)\\
 &\leq&\int_0^t(t-s)^{-a}\,\left(1+  E(\|u_n(s,.)\|_q^p)\right)\, ds\, .
 \end{eqnarray*}
 The rest of the proof, similar to that of the case (ii)(b), is omitted.
\bigskip

{\em Case (i)} Let $u_0\in L^{\infty}(Q)$,  $p\in [1,+\infty[$ and suppose that
the covariance function $f$ satisfies (C3); set  
\[M_n(t)=\sup_{x\in Q} E\big( |u_n(t,x)|^{2p}\big)\, .\]  
We again prove (\ref{0}) and (\ref{dalang}). Since $u_0\in  
L^{\infty}(Q)$, the inequality  
 (\ref{intDH}) proves (\ref{0}). Let $\psi$ be defined by (\ref{defpsi}) and  
 let $\varphi(t)=\psi(t)+\sum_{i=1}^N t^{-\alpha_i+\alpha}\in L^1([0,T])$;  
 then Burkholder's, H\"older's inequalities and (L1) yield  
\begin{eqnarray*}  
M_{n+1}(t,x)&\leq & C_p\, \Bigg[ M_0(t) + E\Big(\Big| \int_0^t  
ds\int_Q dy \int_Q |G(t,x;s,y)|\, |\sigma(s,y,u_n(s,y))|\,  
f(y-z)\\ &&\qquad\qquad \qquad\qquad\qquad \times |G(t,x;s,z)|\,  
|\sigma(s,y,u_n(s,z))|\, dz \Big|^p\Big)\\ &&\qquad + \int_0^t  
\int_Q \sum_{i=1}^N |H_i(t,x;s,y)|\, (1+M_n(s)) dy\,  
 ds\Bigg]\\  
&\leq& C_p\, \left( M_0(t) +  \int_0^t \,  
 \varphi(t-s)\, [1+M_n(s)]\, ds \right)\, .  
\end{eqnarray*}  
This implies (\ref{dalang}) and again  
 Lemma 15 in \cite{D} shows that  
\begin{equation}\label{majPi}  
\sup_n\, \sup_{(t,x)\in [0,T]\times  Q} E(|u_n(t,x)|^{2p}) <+\infty\, .  
\end{equation}  
A similar computation for $\displaystyle \Delta_n(t) = \sup_{x\in  
Q} E\big( |u_{n+1}(t,x)- u_n(t,x)|^{2p}\big)$ and the global  
Lipschitz property (L2) of the coefficients with respect to the  
last variable show that (\ref{delta}) holds. As in   
case 2, usual arguments prove that the solution $u$ to  
(\ref{evolip}) exists and satisfies (\ref{momlip}). 
\epr \\  
\bigskip

\subsection{Cahn-Hilliard equation in dimension $d=4,5$}  
The following stochastic Cahn-Hilliard equation has been studied  
in dimension 1 up to 3 by C.~Cardon-Weber \cite{CW1} and  
\cite{CW2}; see also G. Da Prato and A. Debussche \cite{DPD}. Let  
$Q=[0,\pi]^d$ or a compact convex subset of $\mathbb{R}^d$, 
$(t,x)\in [0,T]\times Q$ and multi-indices $(k_i\, ,  
\, 1\leq i\leq N)$ with $|k_i|\leq 3$ which satisfy
the conditions in Remark \ref{Excond}, the following equation is  
defined in a weak sense:  
\begin{eqnarray}\label{eq:spde1}  
\nonumber \frac{\partial u}{\partial t}(t,x) +(\Delta ^2 u(t,x)  
-\Delta R(u)(t,x) )&=&  
 \sigma(t,x,u(t,x))\dot F+g(t,x,u(t,x))\\&& +\sum_{i=1}^ND^{k_i}_x  
 (b_i(t,x,u(t,x)))\, ,  
 \end{eqnarray}  
with the  initial condition $u(0,.)=u_0$ and the homogeneous Neumann  
boundary conditions:  
 \begin{equation}\label{eq:cond}  
 \frac{\partial  u} {\partial n}= \frac{\partial \Delta u}{\partial n}=0  
 \textrm{ on } \partial Q.  
 \end{equation}  
We will also consider the homogeneous Dirichlet boundary  
conditions:  
 \begin{equation}\label{eq:cond2}  
  u =  \Delta u=0 \textrm{ on } \partial Q.
 \end{equation} 
In this section, we will suppose  
that $d=4,5$ and make the following assumptions:\smallskip

\noindent \textbf{(H.1)} $R$ is a polynomial of degree $3$ with  
positive dominant  coefficient, and such that $R(0)=0$ if (\ref{eq:cond2})
holds.\smallskip

\noindent \textbf{(H.2)}  
  $\sigma:[0,T]\times Q\times\bR \mapsto \bR $ is   bounded,  
  and the functions $\sigma$ and $b_i: [0,T]\times Q\times\bR \mapsto \bR $  
   are globally Lipschitz with respect to the last variable, while  
   the function $g:[0,T]\times Q\times\bR \mapsto \bR $ has  
   quadratic growth with respect to the last variable uniformly  
   with respect to the first ones, and satisfies for $y,z\in Q$:  
   \begin{equation}\label{accroisg} \sup_{(t,x)\in [0,T]\times Q}  
   |g(t,x,y)-g(t,x,z)|\leq C\,  
   (1+|y|+|z|)\, |y-z|\, .\end{equation}  
  
\noindent \textbf{(H.3)}  $u_0$ belongs to $L^q(Q)$ for some $q  
>d$.\smallskip

\noindent \textbf{(H.4)} Set  
 $H_i(t,x;s,y)=D^{k_i}_x G(t,x;s,y)$ for  $|k_i|\leq 3$, where $G$ is the Green function associated  
  with the operator $\frac{\partial}{\partial t} + \Delta^2$  
on $Q$ with the homogeneous boundary conditions (\ref{eq:cond}) or  
(\ref{eq:cond2}), and the multi-indices $k_i$ and the functions $b_i$ satisfy the
assumptions in Example \ref{Excond}. \medskip

 The upper estimates of $G$ stated in condition  
(C2) are given in the introduction: $\alpha =\frac{d}{4}$, $  
\beta=\frac{4}{3}$, $\gamma=\frac{1}{3}$, $\delta=\frac{1}{4}$ and  
$\eta=1$. Clearly for $|k_i|\leq 3$, $\alpha_i\leq \alpha  
+\frac{3}{4}$ and (C4) holds. As explained in the introduction, the Green formula  
shows that the weak formulation of (\ref{eq:spde1}) is equivalent with the evolution  
formulation: for $x \in Q, t \in [0,T]$:  
 \begin{eqnarray}\label{eq:u}  
\nonumber &&u(t,x)\; =\; \int_Q G(t,x;0,y)u_0(y)dy+ \int_0^t\int_Q  
 G(t,x;s,y)\sigma(s,y,u(s,y))F(ds, dy)\\  
 &&\qquad + \int_0^t \int_Q \Delta  
 G(t,x;s,y)R(u(s,y))dyds\\  
 & &\qquad +\int_0^t\int_Q \Big[ \sum_{i=1}^N H_i(t,x;s,y)b_i(s,y,u(s,y))  
 + G(t,x;s,y)g(s,y,u(s,y))\Big]\, dyds\nonumber \, .  
 \end{eqnarray}  
The following theorem completes the existence and uniqueness of  
the solution to (\ref{eq:spde1}) in dimension 4 and 5.  
\begin{Th}.\label{thCH}  
Let $Q$ denote either $[0,\pi]^d$ or a compact subset  of 
$\mathbb{R}^d$ with  boundary of class $\mathcal{C}^{4+\lambda}$ for $\lambda >0$, $d=4,5$  and
assume that (H.1)- (H.4) hold. Let $f$ be the  
covariance function of the Gaussian noise $F$ defined by  
(\ref{cov}), which satisfies (C1) and such that for $\e\in ]0,1[$,  
\begin{equation}\label{covCH}  
\int_{B_d(0,1)} f(v)\, |v|^{-d(1+\e)+4}\,  dv <\infty.  
\end{equation}  
There exists a unique adapted process $u$ in  
$L^\infty([0,T],L^q(Q))$ that satisfies  equation (\ref{eq:u}).  
\end{Th}  
\begin{rk}
Under assumptions (H.1), (H.3) and (H.4), the existence result proved in 
\cite{CW1} in dimension 1-3 
(respectively Theorem \ref{thCH}) extends to a 
compact $Q$ with boundary of class $\mathcal{C}^{4+\lambda}$ for $\lambda >0$,
for the differential operator $a(t)\Delta^2$ where the function $a$ is such that 
$\sup_{0\leq t\leq T}a(t)<0$,
and when $F$ is the space-time white noise (respectively when (H.2) holds). 
\end{rk}
\bpr To prove this theorem we at first prove the existence of a  
 solution when the coefficients are truncated. Let  
$K_n:\bR^+ \rightarrow \bR$ be a ${\cal C}^1$  function such that  
\begin{equation}\label{Kn}  K_n(x)=1 \, \textrm{  if }\,  x<n,  
\quad K_n(x)=0 \, \textrm{ if }\, x\leq n+1,\quad \vert K_n \vert  
\leq 1\; \textrm{ and } \vert K'_n \vert \leq 2\, .  
\end{equation}  
We denote by $u_n$ the solution to the following  evolution  
equation with truncated coefficients:  
\begin{eqnarray}\label{eq:un}  
&&u_n(t,x) =  \int_Q G(t,x;0,y)u_0(y)dy + \int_0^t\int_Q  
G(t,x;s,y)\, \sigma(s,y,u_n(s,y))F(ds,dy)\\ && \qquad +\int_0^t  
\int_Q K_n(\Vert u_n(s,.)\Vert_q))\big[\Delta  
G(t,x;s,y)R(u_n(s,y))+  
 G(t,x;s,y)\, g(s,y,u_n(s,y))\big]\, dyds\nonumber\\  
 && \qquad +\sum_{i=1}^N\int_0^t\int_Q H_i(t,x;s,y)\,  
b_i(s,y,u_n(s,y))dy\, ds \, .\nonumber  
\end{eqnarray}  
 Given an adapted process $u$, let $L(u)$ be defined by  
$$ L(u)(t,x)= \int_0^t\int_Q G(t,x;s,y)\, \sigma(s,y,u(s,y))\,  
F(ds,dy).$$ The arguments used in the proof  of  Theorem  
\ref{existlip} show that  $u_n$ exists, is unique and that  for  
any %$p\in [q,+\infty[$,  
 $ p\in[q, \frac{q}{1-\e}]$,  
$$ \sup_{t\in [0,T]} E\big( \|u_n(t,.)\|_q^p\big) < +\infty.$$  
Indeed, let $\mathcal{K}$ denote the set of adapted  
$L^q(Q)$-valued processes such that  for $q\leq p\leq  
\frac{q}{1-\varepsilon}$,  $\|u\|_{\mathcal{K}}^p=\sup_{0\leq  
t\leq T} E(\|u(t,.)\|_q^p)<+\infty$ . For any $u\in \mathcal{K}$,  
$(t,x)\in [0,T]\times Q$, set  
\begin{eqnarray*}  
H_n(u)(t,x)&=&\int_0^t\int_Q \Delta G(t,x;s,y) \,  
K_n(\|u(s,.)\|_q)\, R(u(s,y)) \,dyds\, ,\\  
J_n(u)(t,x)&=&\int_0^t\int_Q G(t,x;s,y)\,K_n(\|u(s,.)\|_q)\,  
g(s,y,(u(s,y)) dyds\, ,\\ 
B(u)(t,x)&=&\sum_{i=1}^N\int_0^t\int_Q  
H_i(t,x;s,y) \, b_i(s,y,u(s,y))\, dyds\, .  
\end{eqnarray*}  
Since $\sigma$ is bounded and (\ref{covCH}) implies  
(\ref{Acovdiff}), Burkholder's  
 inequality and Lemma \ref{CNS} yield that for any adapted process $u$ and  
 $2\leq p<+\infty$, $\sup\{E(|L(u)(t,x)|^p)\, : (t,x)\in [0,T]\times Q\}  
 <+\infty$, so that $\|L(u)\|_{\mathcal{K}}<\infty$. Furthermore, given $u,v\in  
 \mathcal{K}$,  
 (H.2), the fact that (\ref{covCH}) implies (C'3)(q,p) for $q\leq p\leq  
  \frac{q}{1-\e}$ and the argument used to show (\ref{majoTn1}) in the proof  
   of Theorem~\ref{existlip} yield  
\[ E\left(\|L(u)(s,.)-L(v)(s.)\|^p_q\right)\leq C\left(\int_0^T \psi_p(s)\, ds  
\right)^{\frac{p}{2}}\, \sup_{0\leq s\leq T}  
E(\|u(s,.)-v(s,.)\|_q^p)\, .\] Since $\psi_p$ is integrable, $L$  
is a contraction of $\mathcal{K}$ for small enough $T$.  
 For the polynomial term $H_n$, we just need to notice that  
 if $u$ and $v$ belong to  $L^\infty([0,T],L^q(Q))$  
 \begin{equation}\label{eq:difK}  
\Big \Vert K_n(\Vert u(s,.)\Vert_q)R(u(s,.))-K_n(\Vert  
v(s,.)\Vert_q)R(v(s,.)) \Big\Vert_{\frac{q}{3}}\leq  
 C_n \Vert u(s,.)-v(s,.)\Vert_q.  
 \end{equation}  
Using (\ref{momentqJ(v)}), (\ref{eq:difK}) and H\"older's  
inequality, we obtain that for $d<q\leq p<+\infty$,  
\begin{eqnarray*}\sup_{0\leq t\leq  
T}\|H_n(u)(t,.)-H_n(v)(t,.)\|_{\mathcal{K}}^p&\leq& C_n\,  
T^{(p-1)\left(\frac{1}{2}-\frac{d}{2q}\right)}\int_0^T  
(t-s)^{-\frac{d+2}{4}+\frac{d}{4}(1-\frac{2}{q})} \\ &&\qquad  
\times E(\|u(s,)-v(s,.)\|_q^p)\, ds\\ &\leq&C_n \,  
T^{p\left(\frac{1}{2}-\frac{d}{2q}\right)}\,  
\|u-v\|_{\mathcal{K}}^p\, .  
\end{eqnarray*}  
A similar computation based on the quadratic growth and increments  
property of $g$ with respect to the third variable shows that for  
$\frac{d}{4}<q\leq p<+\infty$,  
\[\sup_{0\leq t\leq  
T}\|J_n(u)(t,.)-J_n(v)(t,.)\|_{\mathcal{K}}\leq C_n\,  
T^{1-\frac{d}{4q}}\, \|u-v\|_{\mathcal{K}}\, .\]  
Finally, the estimation of $T_n^2$ made in the proof of Theorem  
\ref{existlip} shows that  
\[\sup_{0\leq t\leq  
T}\|B(u)(t,.)-B(v)(t,.)\|_{\mathcal{K}}\leq C\, \sup_{1\leq i\leq  
N} T^{1+\alpha-\alpha_i}\, \|u-v\|_{\mathcal{K}}\, .\] Hence,  
there exists $T_0>0$, independent of the initial condition $u_0$,  
such that for $0<T\leq T_0$, $L+H_n+J_n+B$ is a contraction of  
$\mathcal{K}$, and hence admits a unique fixed point such that  
$u(0,.)=u_0$. A concatenation argument implies that (\ref{eq:un})  
has a unique solution on $[0,T]$ for an arbitrary terminal time  
$T$.  
  
To prove the existence and uniqueness of $u$ we follow the proof  
in C. Cardon-Weber \cite{CW1}. Let $\tau_n$ be the stopping time  
defined by: $\tau_n=\inf\{t \geq 0, \Vert u_n(t,.)\Vert_q\geq  
n\}$.  
 By uniqueness of the solution to (\ref{eq:un}), the local property of the stochastic  
integrals yields for $m>n$, $u_m(t,.)=u_n(t,.)$ if $t \leq \tau_n$, so that we can define  
a process $u$ by setting  $u(t,.)=u_n(t,.)$ on $t \leq \tau_n$. Set $\tau_\infty=\lim_n \tau_n$.  
Then $u$ is the unique solution of (\ref{eq:spde1}) on the interval $[0,\tau_\infty)$.  
  We just need to prove that $\tau_\infty=+\infty $ a.s.\\  
 Set   $v_n=u_n-L(u_n)$; then for  
   every $T>0$,  $v_n$ is the weak solution on $[0,T]$  
  to the SPDE (with the same boundary conditions as (\ref{eq:spde1})):  
 \begin{equation}  
\left\{\begin{array}{ll}  
 & \frac{\partial v_n }{\partial t}(t,x) +\Delta^2 v_n(t,x)%\\  
 %& \quad  
 -\Delta\Big[K_n(\Vert  
 v_n(t,.)+L(u_n)(t,.)\Vert_q)R(v_n(t,x)+L(u_n)(t,x))\Big] =\\  
 &\qquad K_n(\Vert  v_n(t,.)+L(u_n)(t,.)\Vert_q)\,  
 g(t,x,v_n(t,x)+L(u_n)(t,x))\\&\qquad +\sum_{i=1}^ND_x^{k_i}  
 b_i(t,x,v_n(t,x)+L(u_n)(t,x)),\\  
 &  v_n(0,.)=u_0(.),\\  
 &  \frac{\partial v_n}{\partial n}=\frac{\partial\Delta v_n}  
 { \partial n}=0 \;{\rm  (resp.}\;  v_n=\Delta v_n=0{\rm )} \; \textrm{ on }  
 \; \partial Q .\label{eq:vn2}  
\end{array}\right.  
 \end{equation}  
Since $\sigma$ is bounded, the Garsia-Rodemich-Ramsay Lemma (cf.  
eg. \cite{Ga}), (\ref{covCH}) and lemma \ref{ISaccrois}  yield  
that for any $p\in [2,+\infty[$,  
\begin{equation}\label{eq:Ga}  
\sup_n E(\Vert L(u_n)\Vert_\infty^{p})<\infty \, .  
\end{equation}  
Since $u_0$ belongs to $L^q(Q)$,  
\begin{equation}\label{eq:Gu0}  
\sup_{t\leq T}\Vert G_tu_0\Vert_q\leq \Vert u_0\Vert_q.  
\end{equation}  
 We need to prove a uniform upper  
estimate for the drift terms $H_n(u_n)$ and $J_n(u_n)$ (the  
estimation of the other drift term $B(u_n)$, which is easier, will  
be omitted and to lighten the notations, we will assume that  
$b_i=0$, $1\leq i\leq N$). Since the function $\Delta G$ has a  
regularizing effect, we first show that $u_n$ belongs to the sets  
$L^a([0,T],L^q(Q))$ for some well-chosen $a$. \\  
 Let us introduce some notations: set $\mathcal{A} =-\Delta$,  let  
 $<.,.>$ denote the usual scalar product in $L^2(Q)$, let $(e_n\,  
 ,\, n\in \mathbb{N}^d)$ be a basis of $L^2(Q)$ made of eigenfunctions of  
 $\mathcal{A}$ (namely $e_n(x)=$\linebreak[3]$\Pi_{i=1}^d e_{n_i}(x_i)$ and either  
 $e_0=\frac{1}{\sqrt{\pi}}$ , $e_n(x)=\sqrt{\frac{2}{\pi}}\,  
 \cos(nx)$ for $n>0$ in the case of the Neumann boundary  
 conditions (\ref{eq:cond}), or $e_n(x)=\sqrt{\frac{2}{\pi}}\,  
 \sin(nx)$ for $n>0$ in the case of the Dirichlet boundary  
 conditions (\ref{eq:cond2}) if $Q=[0,\pi]^d$). The corresponding eigenvalues are  
 $\lambda_n^2$, where $\lambda_n=\sum_{i=1}^d n_i^2$.  
 For $\mu\neq 0$ and $u \in Dom(\mathcal{A}^\mu)$, let  
 $$ \mathcal{A}^\mu u =\sum_{k \in \bN^{d\star}}\lambda_k^\mu  
 <\e_k,u>\e_k\, ;$$  
  $\mathcal{A}^\mu u $ exists for every $u$ such that $\sum_{k \in  
   \bN^{d, \star}} \lambda_k^{2\mu}<e_k,u>^2<\infty$.  In the  
 sequel, for a function $u:[0,T]\times Q \rightarrow \bR$, we will set  (if
 $e_0$ is a constant eigenfunction):
 $$ m(u)(t)=<\e_0,u(t,.)>=\pi^{-\frac{d}{2}}\int_Q u(t,x) dx \textrm{ and  
}\tilde u(t,y)  =u(t,y)-m(u)(t).$$  
Apply $\mathcal{A}^{-1}$ to the equation (\ref{eq:vn}) and take its scalar  
product in $L^2(Q)$ with $\tilde v_n(t,.)$; this leads to  
  \begin{eqnarray}\label{eq:Av}  
 & & \Vert \mathcal{A}^{-1/2}\tilde v_n(t,.)\Vert_2^2-  
\Vert \mathcal{A}^{-1/2}\tilde v_n(0,.)\Vert_q^2 +\int_0^t\Vert  
 \mathcal{A}^{1/2}v_n(s,.)\Vert_2^2\, ds\\  
 & & + \int_0^t K_n(\Vert v_n(s,.)+L(u_n)(s,.)\Vert_q)  \int_Q\Big[  
R\big(v_n(s,x)+L(u_n)(s,x)\big)\tilde v_n(s,x)\nonumber\\  
&&\qquad\qquad\qquad +g\big(s,x,v_n(s,x)+L(u_n(s,x)\big)\,  
\mathcal{A}^{-1}\tilde{v}_n(s,x)\Big]\, dxds=0.\nonumber  
\end{eqnarray}  
This equation is justified because $v_n$  
 belongs to $L^\infty([0,T]\times Q)$. Using the properties of the  
polynomial $R$,  computations made to obtain (2.19)-(2.21) in  
\cite{CW1}, we obtain that for some $b>0$,  
 \begin{eqnarray}  
   \label{eq:Av2}  
   & & \Vert \mathcal{A}^{-1/2}\tilde v_n(t,.)\Vert_2^2 +\int_0^t  
   \Vert \mathcal{A}^{1/2}v_n(s,.)\Vert_2^2\,  ds  
 \nonumber\\  
& & \qquad   +\frac{b}{4} \int_0^tK_n(\Vert  
v_n(s,.)+L(u_n)(s,.)\Vert_q)  
 \Vert v_n(s,.)+L(u_n)(s,.)\Vert_{4}^{4}\, ds \nonumber\\  
 & & \qquad \qquad \leq \int_0^t C(1+ m(u_0)(s)^{4}  
    +\Vert L(v_n)(s,.)\Vert_{4}^{4})\, ds+\Vert \mathcal{A}^{-1/2}\tilde  
 u_0(.)\Vert_2^2.  
 \end{eqnarray}  
Let us find a second ``a priori'' estimate. Denote by $v_n^m$ the  
Galerkin  
 approximation of $v_n$ and let  $P_m$ be the orthogonal projector on  
$\mbox{\rm  Span}\{e_0,...,e_m\}$. For every $\omega$, $v_n^m$ is  
the ``strong`` solution to the following PDE:  
 \begin{equation}\label{eq:vn}  
\left\{\begin{array}{ll}  
 & \frac{\partial v_n^m }{\partial t}(t,x) +\Delta^2 v_n^m(t,x) \\  
& \quad -\Delta\Big[K_n(\Vert  
 v_n^m(t,.)+L(u_n)(t,.)\Vert_q)  
 P_m(R(v_n^m(t,x)+L(u_n)(t,x)))\Big]\\  
 &- K_n(\Vert  
 v_n^m(t,.)+L(u_n)(t,.)\Vert_q)  
 P_m(g(t,x,v_n^m(t,x)+L(u_n)(t,x)))=0\, .  
 \end{array}  \right .  
\end{equation}  
The solution $v_n^m$ to (\ref{eq:vn}) is unique on some random  
time interval $[0,t_n^m[$ and we prove that $t_n^m=+\infty$.  The  
boundary conditions satisfied by $v_n^m$ and the Green Formula  
yield  
 $$\int_Q \Delta^2 v_n^m(t,x)\times v_n^m(t,x)\, dx =\Vert \Delta  
v_n^m(t,x)\Vert_2^2 .$$  
 We now take the  scalar product in  
 $L^2(Q)$ of (\ref{eq:vn}) with $v_n^m$; using once more the Green formula, we  
obtain  
 \begin{eqnarray*}  
 && \frac{1}{2} \frac{\partial }{\partial  
t}\Vert v_n^m(t,.)\Vert_2^2 +\int_Q \Delta^2 v_n^m(t,x)\times  
v_n^m(t,x)\, dx =K_n(\Vert  
 v_n^m(t,.)+L(u_n)(t,.)\Vert_q)\\  
&& \times\int_Q \Big[ R(v_n^m(t,x)+L(u_n)(t,x))) \, \Delta  
v_n^m(t,x)+g(t,x,v_n^m(t,x)+ L(u_n)(t,x))\, v_n^m(t,x)\Big]\,  
dx.%\\  
%&&\qquad\qquad \qquad +g(t,x,v_n^m(t,x)+L(u_n)(t,x))\,  
%v_n^m(t,x)\Big]\, dx=0.  
\end{eqnarray*}  
Using the local Lipschitz property of $R$ and $g(t,x,.)$ and the  
fact that $\int_Q(v_n^m(t,x))^3\Delta v_n^m(t,x)\, dx $ is  
negative, that the  leading coefficient of $R$ is positive and  
that $\Vert K_n\Vert_\infty \leq 1$, we obtain:  
 \begin{eqnarray}\nonumber  
& & \Vert v_n^m(t,.) \Vert_2^2 + \int_0^t \Big[\Vert \Delta  
v_n^m(s,.) \Vert_2^2 +m(v_n^m(s,.))^2\Big]\, ds\leq  
 \Vert u_0\Vert_2^2  
 +C_T\, (1+  \Vert L(u_n)\Vert_\infty^6+ m(u_0)^4)\\  
% \nonumber &&\qquad \leq  \Vert u_0\Vert_2^2  
% +C_T\, (1+  \Vert L(u_n)\Vert_\infty^6+ m(u_0)^2)\\  
 \label{eq:m2} & & \qquad \qquad + C (1+ \Vert L(u_n)\Vert_\infty^2)  
 \int_0^t \Vert v_n^m(s,.) \Vert_4^4K_n(\Vert v_n^m(s,.)+L(u_n)(s,.)\Vert_q)  
 \, ds.  
\end{eqnarray}  
 The norm $(\Vert \Delta\bullet\Vert_{L^2(Q)}^2+ m(\bullet)  
 ^2)^{\frac{1}{2}}$ is equivalent with the Sobolev  
 norm of $W^{2,2}(Q)$ (cf. eg. Da~Prato and Debussche (1996) p.~245).  
The sequence $(v^m_n)_{m }$ is bounded in  
 $L^2([0,T],W^{2,2}(Q))$. Thus,  $t_n^m=\infty$ and this sequence  
 converges as $m\rightarrow +\infty$ in the weak$^\star$ topology of  
 $L^2([0,T],$\linebreak[3]$W^{2,2}(Q))$ . Its weak  
 limit is  the weak solution to (\ref{eq:vn2}) and hence is  
equal to  $v_n$. Therefore $v_n$ belongs to $L^2([0,T],  
W^{2,2}(Q))$, and  
 we can repeat  the  above  computation with $v_n$  instead of  
 $v_n^m$, which yields  
\begin{eqnarray*}  
& & \Vert v_n(t,.) \Vert_2^2 + \int_0^t \Big[\Vert \Delta v_n(s,.)  
\Vert_2^2 +m(v_n(s,.))^2\Big]\, ds  
 \leq  \Vert u_0\Vert_2^2  
 +C_T(1+  \Vert L(u_n)\Vert_\infty^6+ m(u_0)^4)\\  
 & & \qquad + C (1+ \Vert L(u_n)\Vert_\infty^2)  
 \int_0^t \Vert v_n(s,.) \Vert_4^4K_n(\Vert v_n(s,.)+L(u_n)(s,.)\Vert_q)  
 \, ds.  
\end{eqnarray*}  
 Thus, (\ref{eq:Av2}) and  
Schwarz's inequality  imply that  
\begin{eqnarray*}  
& & \Vert v_n(t,.) \Vert_2^2 + \int_0^t \Big[\Vert \Delta v_n(s,.)  
\Vert_2^2 +m(v_n(s,.))^2\Big]\, ds  
 \leq  \Vert u_0\Vert_2^2 \\  
 & & \qquad +C_T(1+  \Vert L(u_n)\Vert_\infty^6)  
 +C_T\,  (1+ \Vert L(u_n)\Vert_\infty^2)\Big[\Vert \mathcal{A}^{-\frac{1}{2}} u_0\Vert_2^2+  
  m(u_0)^4\Big].  
\end{eqnarray*}  
Inequality (\ref{eq:Ga}) yields that for $\beta\in ]1,+\infty[$,  
\begin{eqnarray}  
\label{eq:R1} \sup_n E(\sup_{t \in [0,T]} \Vert v_n(t,.)\Vert_2^{2  
\beta})< \infty, \\  
 \label{eq:R2} \sup_n E\Big(\Big[\int_0^{T} \Big\{\Vert \Delta  
 v_n(t,.)\Vert_2^2 + m(v_n(t,.))^2\Big \}  
 \, dt \Big]^\beta\Big)< \infty.  
\end{eqnarray}  
  
 Moreover, by Sobolev's embedding theorem (Adams 1975, Corollary 5.16) there  
 exists $C>0$ such that for  
 $d\geq 4$ and  $2 \leq r \leq \frac{2d}{d-4}$,  
if $u\in W^{2,2}(Q)$, $ \Vert u\Vert_{L^r(Q)} \leq C \Vert  
u\Vert_{W^{2,2}(Q)}.$  
  Thus, (\ref{eq:R2}) becomes for  $2 \leq r \leq \frac{2d}{d-4}$,  $\beta\in [1,+\infty[$,  
  \begin{equation}\label{eq:R2b}  
\sup_n E\Big(\Big[\int_0^{T} \Vert v_n(t,.)\Vert_r^2  
\,  dt\Big ]^\beta\Big)< \infty.  
  \end{equation}  
  The inequalities (\ref{eq:Ga}), (\ref{eq:R1}) and (\ref{eq:R2b})  
  imply for $2 \leq r < \frac{2d}{d-4}$,  $\beta\in [1,+\infty[$,  
 \begin{eqnarray}  
\label{eq:u1} \sup_n E(\sup_{t \in [0,T]} \Vert u_n(t,.)\Vert_2^{2  
\beta})< \infty,\\  
 \label{eq:u2} \sup_n E\Big(\Big[\int_0^{T} \Vert u_n(t,.)\Vert_r^2  
\,  dt \Big]^\beta\Big)< \infty.  
\end{eqnarray}  
Let us use the interpolation method to prove that $u_n$ belongs  
a.s. to $L^a([0,T],L^{\bar{R}}(Q))$ for $\frac{2d}{d-4}>r \geq \bar{R}\geq 2$,  
$a \geq  1 \vee \frac{2\bar{R}}{r}$.  Set $\bar{R}=(1-\lambda )\, 2+ \lambda\,  
r$ for $\lambda \in [0,1]$; H\"older's inequality implies that  
 $$\int_0^T \Vert u_n(t,.)\Vert_{\bar{R}}^a \, dt  
\leq \int_0^T \Vert u_n(t,.)\Vert_2^{\frac{2a(1-\lambda)}{\bar{R}}}\:  
\Vert u_n(t,.)\Vert_r^{\frac{a r\lambda}{\bar{R}}} \, dt.$$  
Let   $\lambda=\frac{2\bar{R}}{ar} $, we obtain 
$$ \int_0^T \Vert  
u_n(t,.)\Vert_{\bar{R}}^a \, dt \leq \sup_{t \in [0,T]} \Vert  
u_n(t,.)\Vert_2^{\frac{2}{\bar{R}}a(1-\lambda)} \times \int_0^T \Vert  
u_n(t,.)\Vert_r^2\, dt;$$ 
(\ref{eq:u1}) and (\ref{eq:u2}) imply  
that for  $\bar{R} \in [2, \frac{2d}{d-4}[$  and $a \geq 2$,  
\begin{equation}\label{eq:u3}  
\sup_n E\Big(\Big[\int_0^T \Vert u_n(t,.) \Vert_{\bar{R}}^a\,  
dt\Big]^\beta\Big)<\infty.  
\end{equation}  
 Using lemma \ref{convolution} with  
$\rho=\frac{\bar{R}}{3}$, so that  
$\frac{1}{r'}=1+\frac{1}{q}-\frac{3}{\bar{R}}$, we obtain $$ \Vert  
H_n(u_n)(t,.)\Vert_q \leq C\int_0^t(t-s)^{-\frac{d+2}{4}  
  +\frac{d}{4r'}} (\Vert u_n(s,.)\Vert_{\bar{R}}^{3}+1)\, ds. $$  
Let $\gamma, \gamma' \in ]1,+\infty[$ be conjugate exponents, with  
$\gamma$ close enough to one to ensure $\gamma\, (-\frac{d+2}{4}  
+\frac{d}{4r'})>-1$; this is possible (choosing $\bar{R}$ close enough  
to $\frac{2d}{d-4}$) if  
$\frac{1}{2}+\frac{d}{4q}>\frac{3}{8}(d-4)$, i.e.,  
$(\frac{3}{8}-\frac{1}{4q})d<2$. Since $q>d$, this yields $d\leq  
5$ for any $q\in ]d,+\infty[$.  Then H\"older's inequality implies  
 $$  
 \Vert H_n(u_n)(t,.)\Vert_q \leq  
 C\Big[\int_0^t(t-s)^{(-\frac{d+2}{4} +\frac{d}{4r'})\gamma}\, ds  
 \Big]^{\frac{1}{\gamma}} \, \Big[\int_0^t (\Vert  
 u_n(s,.)\Vert_{R}^{3}+1)^{\gamma'}\,  ds\Big]^{\frac{1}{\gamma'}}. $$  
 Using  
 (\ref{eq:u3}), we obtain  
 \begin{equation}\label{eq:Hnu}  
 \sup_n E\Big(\sup_{t \in [0,T]} \Vert H_n(u_n)(t,.)\Vert_q^\beta\Big)<\infty;  
 \end{equation}  
A similar computation (using the quadratic growth of $g$) yields  
for $\rho=\frac{\bar{R}}{2}$, $\frac{1}{r"}=1-\frac{2}{\bar{R}}+\frac{1}{q}$  
and $\gamma$ close enough to one to ensure that $\gamma d  
(-\frac{1}{4}+\frac{1}{4r"})>-1$ (i.e., for $d<8$ and $\bar{R}$ close to  
$\frac{2d}{d-4}$) yields  
\begin{equation}\label{eq:Mnu}  
\sup_n E\Big(\sup_{t\in  
[0,T]}\|J_n(u_n)(t,.)\|_q^\beta\Big)<+\infty\, .  
\end{equation}  
The equations (\ref{eq:Ga})-(\ref{eq:Gu0}), (\ref{eq:Hnu}) and  
(\ref{eq:Mnu})  imply that for $\beta \in [q,+\infty[ $, $d<q $:  
$$ \sup_n E\Big(\sup_{t \in [0,T]} \Vert  
u_n(t,.)\Vert_q^\beta\Big)<\infty.$$  
 We can now conclude that $\tau_\infty=+\infty $ a.s.; indeed, for every  
 $T>0$,  
 % \begin{eqnarray*}  
\[ P(\tau_n\leq T) =  P\Big(\sup_{t \leq T}\Vert u_n(t,.)\Vert_q\geq n\Big )  
%\\  
  \leq  n^{-\beta}\,  E\Big(\sup_{t \leq T}\Vert u_n(t,.)\Vert_q^{\beta}  
 \Big)\, ,\]  
% \end{eqnarray*}  
 so that  $\lim_{n \rightarrow \infty} P(\tau_n \leq T)=0$. Therefore, we  
 can construct the solution to the SPDE (\ref{eq:u}) on any interval  
  $[0,T]$.  
\epr

\section{Regularity of the solution}

The following lemma studies the H\"older regularity of the term involving the initial condition. There are many possible situations, depending on the boundary conditions, whether  
$\int_QG(t,x;0,y)dy=1$ or not, which requires two different arguments. For  
$(t,x)\in ]0,T]\times Q$, set $G_tu_0(x)=\int_Q G(t,x;0,y)\,  
u_0(y)\, dy $ and set $G_0u_0=u_0$.  
\begin{lem}.\label{u0}  Suppose that $Q$ is convex and that  $G$ 
satisfies (\ref{DDH}) with  
$a,b\in  \{0,1\}$.  
  
1) \;(i) Let  $u_0\in L^q(Q)$ for some $q\in  
[1,+\infty[$; then $Gu_0\in \mathcal{C}(]0,T], L^q(Q))$.  
  
\indent\indent (ii) Let  $u_0$ be bounded; then  
for $0<\lambda<1$ and $0<t_0<T$, $Gu_0\in \mathcal{C}^\lambda  
([t_0,T]\times Q)$.  
  
2) Assume furthermore that  
$\int_Q G(t,x;s,y)\, dy=1$ for all $(s,t,x)\in ]0,T]^2\times Q$ with $s<t$.\\  
  
\indent\indent (i) Let $u_0$ be continuous;  
then $Gu_0\in \mathcal{C}([0,T]\times Q)$.\\  
 \indent \indent (ii) Let   
$u_0\in L^q(Q)$; then $Gu_0\in \mathcal{C}([0,T], L^q(Q))$.\\  
 \indent \indent (iii)  
Let $u_0\in \mathcal{C}^\lambda(Q)$ for some $\lambda \in ]0,1[$;  
then for $0\leq s<t\leq T$, ${\displaystyle \sup_{x\in  
Q}|G_tu_0(x)-G_su_0(x)|}$\linebreak[3]$\leq C (t-s)^{  
\frac{\lambda}{4}}$.  
  
3) Let  $Q=[0,M]^d$, 
$u_0\in \mathcal{C}^\lambda(Q)$ for some $\lambda\in ]0,1[$ and  
suppose that for every $1\leq i\leq d$, if $\hat{x}_i=(x_1,  
\cdots, x_{i-1}, x_{i+1}, \cdots, x_d)$, there exists a function  
$\phi_i:[0,T]\times \mathbb{R}\times [0,M]^{d-1}\times  
[0,M]^{d-1}$ such that  
 \[ G(t,x;0,y)=\phi_i(t, x_i+y_i, \hat{x}_i, \hat{y}_i)+\varepsilon_i  
\phi_i(t, x_i-y_i, \hat{x}_i, \hat{y}_i)\, ,\] with  
$\varepsilon_i\in \{-1,1\}$, with $\sup_{(t,x)\in [0,T]\times  
Q}\int_Q |\phi_i(t,x_i+\epsilon y_i, \hat{x}_i,\hat{y}_i)|\, dy  
<+\infty$ for $\epsilon\in \{-1,+1\}$, and suppose that either one  
of the conditions (a) or (b) holds:  
  
\qquad\qquad (a) $\varepsilon_i=+1$ and $\phi_i(t, r+2M,  
\hat{x}_i, \hat{y}_i)= \phi_i(t, r, \hat{x}_i, \hat{y}_i)$ for  
every $r\in \mathbb{R}$.  
  
\qquad\qquad (b) $u_0(x)=0$ for $x\in \partial Q$.  
  
\noindent Then for any $x,x'\in Q$ one has ${\displaystyle  
\sup_{t\in [0,T]}|G_tu_0(x)-G_tu_0(x')|\leq C\, |x-x'|^\lambda}$.  
\end{lem}  
\begin{rk}.\label{remarque}  This proposition can be used for any 
convex compact subset $Q$; in (1), $G$ need not to be a semi-group.   
If $Q=[0,\pi]^d$ and $G$ is the Green function of the operator  
$\frac{\partial}{\partial t}+\Delta^2=0$, for $u_0\in  
\mathcal{C}^\lambda(Q)$, $0<\lambda <1$, the function $ G_.u_0(.)  
\in \mathcal{C}^{\frac{\lambda}{4},\lambda}([0,T]\times Q)$ under  
the homogeneous Neumann boundary conditions (\ref{eq:cond}), while  
under the homogeneous Dirichlet boundary conditions  
(\ref{eq:cond2}), one has $\sup_{0\leq t\leq  
T}|G_tu_0(x)-G_tu_0(x')|\leq C\, |x-x'|^\lambda$, and for  
$0<t_0\leq s<t\leq T$, one has for any $0<\mu<1$, $\sup_{x\in Q}  
|G_tu_0(x)-G_su_0(x)|\leq C(t_0)\, |t-s|^\mu$.  
\end{rk}  
{\it Proof of Lemma \ref{u0}:} 1) (i) Given $0<t_0\leq s<t\leq T$, $0<\lambda<1$,  
$\|G_tu_0-G_su_0\|_q^q\leq C|t-s|^{\lambda q}\,C(M_1+M_2)$, where  
\begin{eqnarray*}  
M_1&=&\int_Qdx \left| \int_Q\left|\int_0^1 (\theta  
t+(1-\theta)s)^{-(\alpha+\eta)}\,  
\exp\left(-c\frac{|x-y|^\beta}{(\theta  
t+(1-\theta)s)^{\gamma}}\right)d\theta\right|^\lambda  
t^{-\alpha(1-\lambda)}\, u_0(y)\, dy\right|^q\\  
 M_2&=&\int_Qdx \left| \int_Q\left|\int_0^1  
(\theta t+(1-\theta)s)^{-(\alpha+\eta)}\, d\theta \right|^\lambda  
s^{-\alpha(1-\lambda)}\exp\left(-c\frac{|x-y|^\beta}{s^{\gamma}}\right)\,  
u_0(y)\, dy\right|^q\, .  
\end{eqnarray*}  
Clearly, $M_1+M_2\leq t_0^{-c}$ for some $c>0$.  
  
\indent (ii) A similar argument for $q=+\infty$ shows that for  
$0<t_0\leq t<t'\leq T$, $x,x'\in Q$, $0<\lambda<1$ and some $c>0$,  
\[|G_tu_0(x)-G_{t'}u_0(x')|\leq C\, t_0^{-c}\,  
\left(|t'-t|^\lambda +|x-x'|^\lambda\right)\, \]  
  
2)  Since $\int_QG(t,x;s,y)\, dy =1$ for every $s<t$,  
$G_tu_0(x)-u_0(x)=\int_QG(t,x;0,y)\, [u_0(y)-u_0(x)]\, dy$ and since  
$G$ is a semi-group, for $t<t'$,  
$G_{t'}u_0(x)-G_tu_0(x)=\int_QG(t,x;0,y)\, dy $\linebreak[3]$
\int_Q  G(t',y;t,z)[u_0(z)-u_0(y)]\, dz$, so that the study of the  
time-regularity is completed by that at 0 .  
  
\indent (i) One has to check that for any $x\in Q$,  
$G_tu_0(x)-u_0(x)$ converges to 0 as $t\rightarrow 0$. The  
argument, based on the continuity of $u_0$ at $x$, is similar to  
the previous one (see e.g. \cite{CW1}, Lemma 2.1).  
  
\indent (ii) Let $u_0 \in L^q(Q)$, let $(u_0^n)_{n\geq 1}$ be a sequence  
of continuous function converging to $u_0$ in $L^q(Q)$. According to  
(i), $(Gu^n_0)$ belong  $C([0,T]\times Q)$ and it suffices to check  
that $\sup_t \Vert G_tu_0\Vert_q \leq C \Vert u_0\Vert_q$. This  
follows from H\"older's inequality and (\ref{intDH}). 
 
\indent (iii) Using the H\"older continuity of $u_0$, one has  
\begin{eqnarray*}  
\sup_{x\in Q}|G_tu_0(x)-u_0(x)|&\leq&\sup_{x\in Q}C\, \int_Q  
t^{-\alpha}\, \exp\left(-c\frac{|x-y|^\beta}{t^\gamma}\right)\,  
|u_0(y)-u_0(x)|\, dy\\ &\leq&\sup_{x\in Q} C\, \int_Q  
t^{-\alpha}\, \exp\left(-c\frac{|x-y|^\beta}{t^\gamma}\right)\,  
\left(\frac{|x-y|}{t^{\frac{\gamma}{\beta}}}\right)^\lambda \,  
t^{\lambda \frac{\gamma}{\beta}}\, dy%\\ &\leq&C\,  
\leq C\, t^{\lambda \frac{\gamma}{\beta}}\, .  
\end{eqnarray*}  
  
\indent (3) The proof of the space regularity under condition  
(a), which is a straightforward extension of that of \cite{BMM}  
Lemma A.2 and \cite{CW1}, Lemma 2.2, is omitted.  
  
We suppose that (b) holds and compare the function $Gu_0$ at  
points $x=(x_1,\hat{x}_1)$ and $x'=(x'_1, \hat{x}_1)$ with  
$x_1<x'_1$; increments of other components are similarly dealt  
with, and provide the required regularity. Obvious changes of  
variables yield $G_tu_0(x)-G_tu_0(x')=\sum_{i=1}^4 D_i(t,x,x')$,  
where if we set $\tilde{Q}=[0,M]^{d-1}$,  
\begin{eqnarray*}  
&&D_1(t,x,x')=\int_0^{M-(x'_1-x_1)}\int_{\tilde{Q}}\phi_1(t,x'_1+y_1,  
\hat{x}_1,\hat{y}_1)\, [u_0(y)-u_0(y_1+(x'_1-x_1), \hat{y}_1)]\,  
d\hat{y}_1\, dy_1\, ,\\  
&&D_2(t,x,x')=\varepsilon_1\int_{x'_1-x_1}^M\int_{\tilde{Q}}\phi_1(t,x'_1-y_1,  
\hat{x}_1,\hat{y}_1)\, [u_0(y)-u_0(y_1-(x'_1-x_1), \hat{y}_1)]\,  
d\hat{y}_1\, dy_1\, ,\\  
&&D_3(t,x,x')=\int_{M-(x'_1-x_1)}^M\int_{\tilde{Q}}\phi_1(t,x'_1+y_1,  
\hat{x}_1,\hat{y}_1)\, [u_0(y)-u_0(M, \hat{y}_1)]\, d\hat{y}_1\,  
dy_1\\ &&\qquad -  
\int_{-(x'_1-x_1)}^0\int_{\tilde{Q}}\phi_1(t,x'_1+y_1,  
\hat{x}_1,\hat{y}_1)\, [u_0(y_1+(x'_1-x_1),\hat{y}_1)-u_0(0,  
\hat{y}_1)]\, d\hat{y}_1\, dy_1\, ,\\ && D_4(t,x,x')=\varepsilon_1  
\,\int_0^{x'_1-x_1}  
\int_{\tilde{Q}}\phi_1(t,x'_1-y_1,\hat{x}_1,\hat{y}_1)\,  
[u_0(y)-u_0(0, \hat{y}_1)]\, d\hat{y}_1\, dy_1\\ &&\qquad  
-\varepsilon_1\int_M^{M+(x'_1-x_1)}\int_{\tilde{Q}}\phi_1(t,x'_1-y_1,  
\hat{x}_1,\hat{y}_1)\, [u_0(y_1-(x'_1-x_1),  
\hat{y}_1)-u_0(M,\hat{y}_1)]\, d\hat{y}_1\, dy_1\, .  
\end{eqnarray*}  
The H\"older regularity of $u_0$ and the  integrability property  
of $\phi_1(t,.,\hat{x}_1,.)$, uniformly with respect to  
$(t,\hat{x}_1)$, conclude the proof. \epr \\  
  
We suppose that $u_0\in \mathcal{C}^a(Q)$ for some $a\in ]0,1[$;  
then $u_0\in L^q(Q)$ for any $q>d$, so that by Theorem \ref{thCH},  
the solution $u$ to (\ref{eq:u}) belongs to $L^\infty([0,T],  
L^q(Q))$ for every $q\in ]d,+\infty[$. Remark \ref{remarque} gives  
the regularity of $Gu_0$ depending on the boundary conditions,  
while Lemma \ref{ISaccrois} gives the regularity of the stochastic  
integral in (\ref{eq:u}). Thus it suffices  
to study the regularity of the drift terms of (\ref{eq:u}) with  
coefficients  which may have polynomial growth.  
\begin{lem}.\label{driftpoly}  
Let $G$ be a (non-necessarily time-homogeneous) semi-group 
satisfying (\ref{DDH}) with  
$\alpha=\frac{\gamma}{\beta}d$, let $a$ be a multi-index such that  
$|a|\delta<1$, $H(t,x;s,y)=D_x^aG(t,x;s,y)$. Let $b:[0,T]\times  
Q\times \mathbb{R}\rightarrow \mathbb{R}$ be measurable such that  
$\sup\{|b(t,x,y)|\, , \, (t,x)\in [0,T]\times Q\}\leq C\, |y|^m$  
for some $m\geq 1$, let $0<\lambda  
<\frac{1-|a|\delta}{\eta\vee 1}$ and $0<\mu<\left[ (1-|a|\delta)\,  
\big( \frac{\beta}{\gamma}\wedge \frac{1}{\delta}\big) \right ]\wedge 1  
$. Then if  $u:\Omega\times [0,T]\times Q\rightarrow  
\mathbb{R}$  is  a process in $L^\infty([0,T],L^q(Q))$ for  
$q$ large enough,  
  
\[B(u)(t,x)=\int_0^t\int_QH(t,x;s,y)\, b(s,y,u(s,y))\, dy\, ds\,  .\]  
 the map $ B(u)$ belongs to  
$\mathcal{C}^{\lambda,\mu}([0,T]\times Q)$.  
\end{lem}  
\bpr The argument, based on the factorization method (see e.g. G.  
Da Prato and J. Zabczyk \cite{DPZ}) is similar to that in the  
proof of section 2.3 in \cite{CW1}; it is briefly sketched.  
  
Let $\bar{\delta}=|a|\, \delta$,  $\varepsilon\in ]0,1[$ and set  
\begin{eqnarray*}  
\mathcal{J}(v)(t,x)&=&\int_0^t\int_Q G(t,x;s,y)\,  
(t-s)^{-\varepsilon}\, v(s,y)\, dy ds\, ,\\  
\mathcal{K}(v)(t,x)&=&\int_0^t\int_Q H(t,x;s,y)\,  
(t-s)^{\varepsilon-1}\, b(s,y,v(s,y))\, dy ds\, .  
\end{eqnarray*}  
The semi-group property of $G$ implies that for every $(t,x)\in  
[0,T]\times Q$, $B(u)(t,x)=\frac{\sin(\varepsilon  
\pi)}{\pi}$\linebreak[3]$\times \mathcal{J}(\mathcal{K}(u))(t,x)$.  
We prove that for $\varepsilon  
>\bar{\delta}+\alpha\, \frac{m-1}{q}$, $\mathcal{K}$ maps  
$L^\infty([0,T],L^q(Q))$ into itself. Indeed, it suffices to use  
Lemma \ref{convolution} with $q$ and $\frac{q}{m}$, so that  
$\frac{1}{r}=1-\frac{m-1}{q}$, and (\ref{intDH}). We then prove  
that for $v\in L^\infty([0,T],L^q(Q))$, the trajectories of  
$\mathcal{J}(v)$ have the required H\"older regularity. Let  
$x,x'\in Q$; then $|\mathcal{J}(v)(t,x)-\mathcal{J}(v)(t,x')| \leq  
A_1(t,x,x')+A_2(t,x,x')$, where  
\begin{eqnarray*}  
A_1(t,x,x')&=&\int_0^t\int_Q 1_{\{|y-x|\leq |x'-x|\}}  
(t-s)^{-\varepsilon}\, \Big(|G(t,x;s,y)|+|G(t,x';s,y)|\Big)\,  
|v(s,y)|\, dy\, ds\, ,\\ 
A_2(t,x,x')&=&\int_0^t\int_Q 1_{\{|y-x|>  
|x'-x|\}} (t-s)^{-\varepsilon}\, |G(t,x;s,y)-G(t,x';s,y)|\,  
|v(s,y)|\, dy\, ds\, .  
\end{eqnarray*}  
H\"older's inequality and a change of variables yield that for  
$0<\mu<\frac{q-1}{q}$,  
\begin{eqnarray}\label{A1x}  
A_1(t,x,x')&\leq & C\, \int_0^t (t-s)^{-(\varepsilon +\alpha)}\,  
\|v(s,.)\|_q\left(\int_{|z|\leq |x-x'|\,  
t^{-\frac{\gamma}{\beta}}} \exp(-c|z|^\beta)\, t^\alpha\,  
dz\right)^{\frac{q-1}{q}}\, ds\nonumber\\ &\leq&|x-x'|^\mu \,  
\int_0^t (t-s)^{-\varepsilon  
-\frac{\alpha}{q}-\frac{\gamma\mu}{\beta}}\, ds\,  .  
\end{eqnarray}  
The convergence of this last integral requires  
$\mu<\frac{\beta}{\gamma}\left(1-\varepsilon  
-\frac{\alpha}{q}\right)$ and $1-\varepsilon-\frac{\alpha}{q}>0$. Furthermore, if $|x-y|> |x'-x|$, and  
$\tilde{x}$ is a convex combination of $x$ and $x'$, then  
$|\tilde{x}-y|\geq \frac{1}{\sqrt{2}}\big(|2x-x'-y|\wedge |x'-y|\big)$.  
Therefore, Taylor's formula and H\"older's inequality imply that  
for $0<\mu<1$,  
\begin{equation}\label{A2x}  
A_2(t,x,x')\leq C\, |x-x'|^\mu  
\int_0^t (t-s)^{-\varepsilon-\mu \delta - \frac{\alpha}{q}}
\, \|v(s,.)\|_q\, ds\, ,  
\end{equation}  
and the last integral converges if $\mu<\frac{1-\varepsilon  
-\frac{\alpha}{q}}{\delta}$.  
  
Similarly, for $0<t<t'\leq T$ and $x\in Q$,  
$|\mathcal{J}(v)(t,x)-\mathcal{J}(v)(t',x)|\leq  
B_1(t,t',x)+B_2(t,t',x)$, where  
\begin{eqnarray*}  
B_1(t,t',x)&=&\int_0^t\int_Q|(t'-s)^{-\varepsilon}\,  
G(t',x;s,y)-(t-s)^{-\varepsilon}\, G(t,x;s,y)|\, |v(s,y)|\, dy\,  
ds\, ,  
\\ B_2(t,t',x)&=&\int_t^{t'}\int_Q(t'-s)^{-\varepsilon}\,  
|G(t',x;s,y)|\, |v(s,y)|\, dy\, ds\, .  
\end{eqnarray*}  
Computations similar to the previous ones yield for $\lambda \in  
]0,1[$  
\begin{eqnarray}  
\label{B1t} B_1(t,t',x)&\leq&C\, |t-t'|^\lambda\int_0^t  
(t-s)^{-\varepsilon-\eta\lambda -\frac{\alpha}{q}}\,  
\|v(s,.)\|_q\, ds\, ,\\\label{B2t}  
 B_2(t,t',x)&\leq& C\, \int_t^{t'} (t'-s)^{-\varepsilon-\frac{\alpha  
}{q}}\|v(s,.)\|_q\, ds \, \leq C \, \vert  
t-t'\vert^{1-\varepsilon-\frac{\alpha}{q}} .  
\end{eqnarray}

The integral in the right hand-side of (\ref{B1t}) converges if  
$\lambda <\frac{1-\varepsilon-\frac{\alpha}{q}}{\eta}$, while  
that in the right hand-side of (\ref{B2t}) converges for  
$\varepsilon<1-\frac{\alpha}{q}$. Thus for $q$ arbitrary large  
and $\varepsilon$ close enough to $\bar{\delta}$, we see that the  
inequalities (\ref{A1x})-(\ref{B2t}) conclude the proof. \epr  
\bigskip

The following theorem summarizes the results proved in  
Lemmas \ref{ISaccrois}, \ref{driftpoly} and  
Remark \ref{remarque}.  
\begin{Th}.\label{holdCH}  Assume (H.1), (H.2) and (H.4), let $d=4,5$,  $u_0\in \mathcal{C}^a(Q)$ for  
some $a\in ]0,1[$, $F$ be a Gaussian process with covariance  
define in terms of the function $f$ by (\ref{cov}), such that (C1)  
holds and $\int_{B_d(0,1)}f(v)\, |v|^{-d(1+\varepsilon)+4}\,  
dv <+\infty$ for some $\varepsilon>0$. Then the solution  
$u$ to (\ref{eq:u})  
 belongs to $\mathcal{C}^{\bar{\lambda},\bar{\mu}}([0,T]\times Q)$  
 under the Neumann boundary conditions (\ref{eq:cond}) (resp. the  
 map $x\rightarrow u(t,x)\in \mathcal{C}^{\bar{\mu}}(Q)$  
 uniformly in $t\in [0,T]$ while the map $t\rightarrow u(t,x)\in  
 \mathcal{C}^\lambda([t_0,T])$ uniformly in $x\in Q$ for $0<t_0<T$ under  
 the Dirichlet boundary conditions (\ref{eq:cond2})), where  
 \[0<\lambda<\left(1- \frac{\max_ik_i\vee 2}{4}\right)\wedge 
 \frac{\varepsilon d }{8},
 \quad  
 \bar{\lambda}=\lambda\wedge \frac{a}{4} \quad and \quad  
 0<\bar{\mu}<a\wedge \frac{\varepsilon d }{2}. \]
 \end{Th}  
 Finally, a straightforward extension of the preceding computations,
 using Lemmas \ref{ISaccrois} and \ref{driftpoly} (with $\alpha_i$ instead of
 $\alpha +|a|\delta$), provides H\"older regularity for the solution
 to (\ref{evolip}).
 \begin{Th}\label{holdlip} Let $Q$ be convex, suppose that $G$ satisfies
 (\ref{DDH})  and  that
 the assumptions of Theorem \ref{existlip} are fulfilled. Let $u_0$ be
 a continuous function on $Q$, let $u$ be the solution to (\ref{evolip})
 and set $v(t,x)=u(t,x)-G_tu_0(x)$. Then
 
 (i) If $\alpha <1$ and $F$ is the space-time white noise, then $v\in 
 \mathcal{C}^{\lambda,\mu}([0,T]\times Q)$ for $0<\lambda< 
 \inf_i(1+\alpha-\alpha_i) \wedge 
 \frac{1-\alpha}{2}$ and $0<\mu< \frac{\inf_i(1+\alpha-\alpha_i)}{\delta}
\wedge  \frac{1-\alpha}{2\delta}\wedge 1$.
 
 (ii) If $\alpha \geq 1$ and $F$ is a Gaussian process with covariance function
 $f$ defined by (\ref{cov}) such that (C1) holds and $\int_{B_d(0,1)} f(v)\, 
 |v|^{-d(1+\epsilon)+\frac{1}{\delta}}\, dv <+\infty$
 for some $\epsilon >0$, then $v\in \mathcal{C}^{\lambda,\mu}([0,T]\times Q)$ 
 for $0<\lambda< \inf_i(1+\alpha-\alpha_i)
 \wedge (\epsilon d\delta)\wedge 1$
 and $0<\mu< \frac{\inf_i(1+\alpha-\alpha_i)}{\delta} 
 \wedge (\epsilon d)\wedge 1$.
 \end{Th}

\section{Density of the solution to the stochastic Cahn-Hilliard  
PDE}\label{densite} In this section, we concentrate on the  
solution to (\ref{eq:u}) in dimension 4 and 5 under either the  
homogeneous Neumann or Dirichlet boundary conditions on $Q=[0,\pi]^d$. Thus,
 we prove that under proper non-degeneracy conditions on  
the "diffusion" coefficient $\sigma$, the law of $u(t,x)$ has a  
density for $t>0$ and $x\in Q$. This extends results  
proved in \cite{CW1} and \cite{CW2} to higher dimension. Since the  
noise $F$ has a space correlation, the setting of the  
corresponding stochastic calculus of variations is that used in  
\cite{MS}. Let $Q=[0,\pi]^d$,  $\mathcal{E}$ denote the inner  
product space of measurable functions $\varphi:Q\rightarrow  
\mathbb{R}$ such that $\int_Q dx\int_Q dy\,  |\varphi(x)|\:  
f(x-y)\: |\varphi(y)|\, <+\infty$,  endowed with the inner product  
\[<\varphi,\psi>_{\mathcal{E}} =\int_Q dx\int_Q dy \, \varphi(x)\,  
f(x-y)\, \psi(y)\, .\] Let $\mathcal{H}$ denote the completion of  
$\mathcal{E}$ and set $\mathcal{E}_T=L^2([0,T], \mathcal{E})$ and  
$\mathcal{H}_T=L^2([0,T], \mathcal{H})$. Note that $\mathcal{H}$  
and $\mathcal{H}_T$ need not be spaces of functions, and that  
$\mathcal{H}_T$ is a Hilbert space which is isomorphic to the  
reproducing kernel Hilbert space of the Gaussian noise  
$(F(\varphi)\, ;\, \varphi\in \mathcal{D}([0,T]\times Q)$. This  
noise can be identified with a Gaussian process $(W(h)\, ,\, h\in  
\mathcal{H}_T)$ defined as follows. Let $(e_j\, , \, j\geq  
0)\subset \mathcal{E}$ be a CONS of $\mathcal{H}$; then  
$(W_j(t)=\int_0^t\int_Q e_j(x)\, F(ds,dx)\, , \, j\geq 0)$ is a  
sequence of independent standard Brownian motions such that  
\[ F(\varphi)=\sum_{j\geq 0} \int_0^t <\varphi(s,*),e_j>_{\mathcal{H}}\,  
dW_j(s)\, , \varphi\in \mathcal{D}([0,T]\times Q)\, .\] For  
$h\in\mathcal{H}_T$, we set $W(h)=\sum_j\int_0^T  
<h(s),e_j>_{\mathcal{H}}\, dW_j(s)$, and use the framework of the  
Malliavin calculus described in \cite{N} to define the Malliavin  
derivative $DX$ of a random variable $X$  and the corresponding  
Sobolev spaces $\mathbb{D}^{N,p}$. Given $X\in \mathbb{D}^{1,2}$,  
$h\in \mathcal{H}_T$, set $D_hX=<DX,h>_{\mathcal{H}_T}= \int_0^T  
<D_{r,*}X,h(r)>_{\mathcal{H}} \, dr$, where $D_{r,*}X\in  
\mathcal{H}$ for every $r\in [0,T]$. Finally for $r\in [0,T]$ and  
$\varphi\in \mathcal{H}$, set  
$D_{r,\varphi}X=<D_{r,*},\varphi>_{\mathcal{H}}$.  
  
Since the coefficients $R$ and $g(t,x,.)$ are locally Lipschitz,  
we need to localize the Sobolev spaces as follows. A random  
variable $X$ belongs to $\mathbb{D}^{1,p}_{\rm loc}$ if there  
exists an increasing sequence $\Omega_n\subset \Omega$ such that  
$\lim_n P(\Omega_n)=1$ and for every $n$, there exists a random  
variable $X_n\in \mathbb{D}^{1,p}$ and $X=X_n$ on $\Omega_n$.  Let  
$u_0\in \mathcal{C}(Q)$ and suppose that the conditions (H1) and  
(H'2) hold, where  
  
\textbf{(H'2)} The function $\sigma:\mathbb{R}\rightarrow  
\mathbb{R}$ is bounded, globally Lipschitz,  the map $ g(t,x,.)$  
is of class $\mathcal{C}^1$ with quadratic growth and satisfies  
(\ref{accroisg}), and the maps $b_i(t,x,.), 1\leq i\leq N$ are of  
class $\mathcal{C}^1$ with derivatives bounded uniformly in  
$(t,x)$.  
  
Let $u$ denote the solution to (\ref{eq:u})  with either the  
homogeneous Neumann or Dirichlet boundary conditions. Lemmas  
\ref{ISaccrois} and \ref{driftpoly}  
imply that the trajectories of $u-Gu_0$ are almost surely H\"older  
continuous on $[0,T]\times Q$, while the function $G u_0$ is  
clearly bounded by $\|u_0\|_{\infty}$. Therefore, $\lim_n  
P(\Omega_n)=1$ if for every $n\geq 1$ one sets  
\[\Omega_n=\Big\{\omega\in \Omega\, :\, \sup\{|u(t,x)|\, , \,  
(t,x)\in Q\}\leq n\Big\}\, .\]  
  
We now construct a sequence of processes $u(n) \in  
\mathbb{D}^{1,p}$ for every $p\in [2,+\infty[$ such that $u=u(n)$  
on $\Omega_n$. Let $K_n$ be the sequence defined at the beginning  
of the proof of Theorem  \ref{thCH}, which satisfies (\ref{Kn});  
set $\bar{R}_n(x)=K_n(|x|) R(x)$ and $\bar{g}_n(t,x,y)= K_n(|y|)\,  
g(t,x,y)$. The functions $\bar{R}_n$ and $y\rightarrow  
\bar{g}_n(t,x,y)$ are of class $\mathcal{C}^1$ with bounded  
derivatives. Hence Theorem \ref{existlip} yields the existence and  
uniqueness of the process $u(n)$ solution to the evolution  
equation:  
\begin{eqnarray}\label{u(n)} u(n)(t,x)&=&G_tu_0(x)+\int_0^t \int_Q  
G(t,x;s,y)\,\sigma(u(n)(s,y))\, F(ds,dy)\\ \nonumber &&+  
\int_0^t\int_Q\Big[ \Delta G(t,x;s,y)\,  
\bar{R}_n(u(n)(s,y))+G(t,x;s,y)\, \bar{g}_n(s,y,u(n)(s,y))\\  
\nonumber &&\qquad \qquad \qquad +\sum_{i=1}^N H_i(t,x;s,y)\,  
b_i(s,y,u(n)(s,y))\Big]\, dy ds\, .  
\end{eqnarray}  
The local property of stochastic integrals implies that $u(n)=u$  
on $\Omega_n$. The following proposition shows that $u(n)\in  
\mathbb{D}^{1,p}$ for every $p\in [2,+\infty[$.  
\begin{pro}.\label{Dlip} Let $Q$ be a compact subset  of  
$\mathbb{R}^d$, $u_0\in \mathcal{C}(Q)$, let $\sigma  
:\mathbb{R}\rightarrow \mathbb{R}$ be globally Lipschitz and for  
$i\in \{1, \cdots , N\}$, let $b_i: [0,T]\times Q\times  
\mathbb{R}$ satisfy the conditions (L1) and (L2), and suppose that  
the maps $y\rightarrow b_i(t,x,y)$ are of class $\mathcal{C}^1$  
with bounded derivatives, let $G$ and $(H_i\, ,\, 1\leq i\leq N)$  
satisfy the conditions (C2) and (C4). Let $F$ be a Gaussian  
process with space covariance defined by (\ref{cov}) in terms of a  
function $f$ which satisfies the condition (C3). Then for every  
$p\in [2,+\infty[$ and $(t,x)\in [0,T]\times Q$, the solution  
$v(t,x)$ to (\ref{evolip}) belongs to $\mathbb{D}^{1,p}$.  
Furthermore, for every $r\in [0,T]$ and $\varphi\in \mathcal{H}$,  
$D_{r,\varphi}v(t,x)=0$ if $r>t$, while there exists a bounded,  
adapted  family of random variables $(S(s,y)\, , \, (s,y)\in  
[0,T]\times Q)$ such that for $0\leq r\leq t$:  
\begin{eqnarray}  
&&D_{r,\varphi}v(t,x)\; =\; \;<G(t,x;r,*)\,  
\sigma(v(r,*)),\varphi>_{\mathcal{H}}+ \int_r^t \int_Q\Big[  
G(t,x;s,y)\, S(s,y)\label{Dv}\\ &&\qquad \qquad \times  
D_{r,\varphi}v(s,y)\, F(ds,dy)+\sum_{i=1}^N H_i(t,x;s,y)\,  
\partial_3 b_i(s,y,v(s,y)) \, D_{r,\varphi}v(s,y)\, dyds\Big]\, ,\nonumber  
\end{eqnarray}  
and for every $p\in [1,+\infty[$,  
\begin{equation} \label{majDv} \sup_{(t,x)\in [0,T]\times Q}  
E\left( \left| \int_0^t \|D_{r,*} v(t,x)\|^2_{\mathcal{H}}\,  
dr\right|^p\right)=C(p)<+\infty\, .  
\end{equation}  
Furthermore, given  
$0\leq s<t\leq T$, if $\psi$ denotes the function defined by  
(\ref{defpsi}),  
\begin{equation}\label{majoDvst} %%????? AMELIORER  
\sup_{x \in  Q} E\left( \int_s^t \|D_{r,\varphi}  
v(t,x)\|_{\mathcal{H}}^2\, dr \right)\leq C \left[  
\int_0^{t-s}\psi(\tau)\, d\tau+  
\sum_{i=1}^N(t-s)^{2(1+\alpha-\alpha_i)}\right]\,  .  
\end{equation}  
\end{pro}  
\bpr Let $(v_k\, ,\, k\geq 0)$ be the Picard approximation scheme  
of $v$ defined by (\ref{picardlip}); then by the proof of Theorem  
\ref{existlip}, the sequence $(v_k(t,x)\, , \, k\geq 0)$ is  
bounded in $L^p(\Omega)$ uniformly in $(t,x)$ and  converges in  
$L^p(\Omega)$ to $v(t,x)$. Following a classical argument, we  
prove by induction on $k$ that $v_k(t,x)\in \mathbb{D}^{1,p}$ and  
that \begin{eqnarray}  
% \begin{equation}  
&& \label{bornDvk} \sup_k \; \sup_{(t,x)\in  
 [0,T]\times Q} E\left( \| D v_k(t,x)\|_{\mathcal{H}_T}^p \right)  
 <+\infty\, ,\\  
% \end{equation}  
% \begin{equation}  
&&\label{serieDvk} \sum_k \sup_{(t,x)\in  
 [0,T]\times Q} E\left(\| Dv_{k+1}(t,x)-Dv_k(t,x)\|^p_{\mathcal{H}_T} \right)  
 <+\infty\, .  
% \end{equation}  
\end{eqnarray}  
 Then using \cite{N}, Lemma 1.2.3, we conclude that $Dv_k(t,x)$  
 converges to $Dv(t,x)$ in the weak topology of $L^p(\Omega,  
 \mathcal{H}_T)$; furthermore, this yields (\ref{majDv}).  
 Since $v_0$ is deterministic, it belongs to $\mathbb{D}^{1,p}$; suppose  
that $v_k\in \mathbb{D}^{1,p}$; since $\sigma$ is globally  
Lipschitz, Proposition 1.2.3 in \cite{N} implies that  
$\sigma(v_k(s,y))\in \mathbb{D}^{1,p}$ and that  
$D_{r,\varphi}(\sigma(v_k(s,y)))= S_k(s,y)\, D_{r,\varphi}  
v_k(s,y)$, where $S_k(s,y)$ is a bounded adapted process.  
Furthermore, for every $r\in [0,T]$ and $\varphi \in \mathcal{H}$,  
$D_{r,\varphi} v_{k+1}(t,x)=0$ if $r>t$ and for $r\leq t$:  
\begin{eqnarray*}  
&&D_{r,\varphi}v_{k+1}(t,x)=<G(t,x;r,*)\, \sigma(v_k(r,*)),\varphi  
>_{\mathcal{H}}+ \int_r^t \int_Q \Big[ G(t,x;s,y)\, S_k(s,y)\,\\&&  
\quad \times D_{r,\varphi}v_k(s,y)\, F(ds,dy)+\sum_{i=1}^N  
H_i(t,x;s,y)\,  
\partial_3 b_i(s,y,v_k(s,y)) \, D_{r,\varphi}v(s,y)\, dyds\Big]\, .  
\end{eqnarray*}  
Let $\psi$ be the $L^1([0,T])$ function defined by (\ref{defpsi})  
and set $I= \int_0^T \psi(r)\, dr$ . The linear growth condition  
on $\sigma$ and equations (\ref{majPi}) and (\ref{inter}) imply  
that for any $p\in [ 2,+\infty[$, there exists a constant $C_p$  
(which does not depend on $k$) such that for every $k$  
\begin{equation}\label{sigmaDvk} \sup_{x\in Q} E\left(  
 \|G(t, x;.,*)\,  
\sigma(v_k(.,*))\|_{\mathcal{H}_T}^{2p}\right) \leq C\, I^{p-1}  
\int_0^t \psi(t-r)\, \sup_{(r,x)} E(|v_k(r,x)|^{2p})\, dr \leq C_p\,  
.  
\end{equation}  
For $t \in [0,T]$, $x\in Q$, let $(Y_{\tau}(x)\, , \, 0\leq \tau \leq T)$ 
be the   $\mathcal{H}_T$-valued martingale defined by  
\[ Y_\tau(x) = \int_0^{\tau\wedge t} \int_Q G(t,x;s,y)\, 
S_k(s,y)\, Dv_k(s,y)\,  
F( ds,dy)\, .\] Let $(\epsilon_j\, , \, j\geq 0)$ be a CONS of  
$\mathcal{H}_T$; then Burkholder's inequality for Hilbert-valued  
martingales (see e.g. \cite{Met}, p. 212) and Parseval's identity  
yield  
\begin{equation}\label{ISDvk} \sup_{x\in Q}  
\|Y_t(x)\|^{2p}_{L^{2p}(\Omega, \mathcal{H}_T)} \leq C_p\,  
I^{p-1}\, \int_0^t \psi(t-s)\, \sup_{y} E(\|D  
v_k(s,y)\|^{2p}_{\mathcal{H}_T})\, ds\, .  
\end{equation}  
Finally, Lemma \ref{convolution}  and the conditions (C2) - (C4)  
imply that the function $\bar{\psi}(s)=\psi(s)+\sum_{i=1}^N  
s^{-\alpha_i +\alpha}\in L^1([0,T])$, and together with the  
inequalities (\ref{sigmaDvk}) and (\ref{ISDvk}) this yields that  
there exists $C_p>0$ such that for every $k\geq 0$,  
\[\sup_{x\in Q}\|Dv_k(t,x)\|_{L^{2p}(\Omega,\mathcal{H}_T)}^{2p}  
\leq C_p + C_p \int_0^t\bar{\psi}(t-s)\, \sup_{y\in Q}  
\|Dv_k(s,y)\|^{2p}_{L^{2p}(\Omega, \mathcal{H}_T)}\, ds\, .\] Thus  
Lemma 15 in \cite{D} concludes the proof of (\ref{bornDvk}). A  
similar argument shows (\ref{serieDvk}). We conclude that each  
$v_k(t,x) \in \mathbb{D}^{1,p}$, and that (\ref{Dv})and  
(\ref{majDv}) hold. To prove (\ref{majoDvst}), we use  
(\ref{momlip}), (\ref{Dv}) and arguments similar to the previous  
ones; then for $0\leq s<t\leq T$,  
\begin{eqnarray*}  
 E \left(  \int_s^t \|G(t,x;r,*)\,  
\sigma(v(r,*))\|^2_{\mathcal{H}}\, dr \right) &\leq& C\, \left(  
\int_s^t \| |G(t,x;r,*)| \|_{\mathcal{H}}^2\, dr \right)\,  
\left[ 1+ \sup_{(r,y)} E(|v(r,y)|^{2p})\right]\\ 
& \leq & C\, \int_0^{t-s}   \psi(\tau )\, d\tau \, .  
\end{eqnarray*}  
Furthermore, the isometry of Hilbert-spaced valued martingales,  
Fubini's theorem and (\ref{majDv}) imply  
\begin{eqnarray*}  
&& E\left(  \int_s^t \left\|\int_r^t \int_Q G(t,x;\tau ,y)\,  
S(\tau,y)\, D_{r,*}v(\tau ,y)\, F(d\tau ,dy)  
\right\|_{\mathcal{H}}^2\, dr\right)\\&&\qquad \leq C\, \int_s^t  
d\tau \int_Q dy \int_Q dz\,  G(t,x;\tau ,y)\,f(y-z)\,  
G(t,x;\tau ,z)\,  
\\&&\qquad\qquad \times \int_s^{\tau}  E(<D_{r,*} v(s,y)\, , \, D_{r,*}  
v(s,z>_{\mathcal{H}})\,dr \leq C\, \int_0^{t-s} \psi(\tau)\,  
d\tau\, .  
\end{eqnarray*}  
Finally,  conditions (C2) and (C4), Minkowski's and Schwarz's  
inequalities and Fubini's theorem imply that for every $i\leq N$,  
\begin{eqnarray*}  
&&E\left(  \int_s^t \left\|\int_r^t \int_Q H_i(t,x;\tau ,y)\,  
\partial_3 b_i(\tau,y,v(s,y)) \, D_{r,*}v(\tau ,y)\, dy d\tau  
 \right\|_{\mathcal{H}}^2\, dr\right)\\  
 &&\qquad \leq C\, (t-s)^{\alpha +1-\alpha_i} \int_s^t d\tau  
 \int_Q dy \, |H_i(t,x;\tau ,y)|\, \sup_{(\tau,y)\in [0,T]\times Q}  
 \int_s^{\tau} E(\| D_{r,*}v(\tau,y)\|_{\mathcal{H}}^2)\, dr\\&&\qquad\leq  
 (t-r)^{2(\alpha +1-\alpha_i)}\, .  
 \end{eqnarray*}  
 This completes the proof of (\ref{majoDvst}).\epr  
  
The following theorem, which establishes the absolute continuity  
of the law to the stochastic Cahn-Hilliard PDE, is the main result  
of this section.  
\begin{Th}.\label{densiteCH} Let $Q=[0,\pi]^d$, suppose that the conditions 
(H.1),   (H'2) and (H.4)  hold, let $F$ be a Gaussian noise with covariance defined by   (\ref{cov}) in terms of a function $f$ which satisfies the  
conditions (C1) and (\ref{covCH}). Let $u_0\in \mathcal{C}(Q)$  
and let $u$ be the solution to (\ref{eq:u}) with initial  
conditions $u_0$ and  the homogeneous Neumann or Dirichlet 
boundary conditions. Let $t_0\in ]0,T]$, $x_1, \cdots, x_l$ be  
pairwise distinct points in $]0,\pi[^d$, and set  
$u(t_0,\underline{x})= (u(t_0,x_1), \cdots ,u(t_0, x_l))$. For any  
$\tau >0$, let  
\begin{equation}\label{Ie}  
I(\tau)=\int_{B_d(0,\tau)} f(v)\, |v|^{-d+4}\,  
\ln\big(|v|^{-1}\big)^{(5-d)^+}\, dv\, .  
\end{equation}  
  
(i) Suppose that $|\sigma|\geq C >0$; then if for some  
$0<\nu<\frac{1}{4}$  
  \[ \lim_{\tau \rightarrow 0}  
I(\tau^{\frac{1}{4}+\nu})^{-1}\, \Big[ \tau + \tau^{\frac{1}{2}}\,  
I(\tau^{\frac{1}{4}-\nu}) \Big] =0\, ,\]  
 the law of  
$u(t_0,\underline{x})$ is absolutely continuous with respect to  
Lebesgue's measure.  
  
(ii) Let $u_0\in \mathcal{C}^a(Q)$ for some $a>0$ and suppose that  
for some $\nu>0$:  
\begin{equation}\label{Ielip}  
\lim_{\tau \rightarrow 0} I(\tau^{\frac{1}{4}+\nu})^{-1}\, \Big[  
\tau^{A}+\tau^{\frac{1}{2}}\,  
I(\tau^{\frac{1}{4}-\nu})\Big] =0\; for\; A\in  
]0,\frac{d\varepsilon}{4}+\frac{\mu}{8}[\; and\; \mu\in ]0,1\wedge  
\frac{d\varepsilon}{2}\wedge a[\, .  
\end{equation}  
Then the law of $u(t_0,\underline{x})$ is absolutely continuous  
%with respect to the Lebesgue measure  
on $\{\sigma\neq 0\}^l$.  
\end{Th}  
\begin{rk}.\label{exemple} Let $f(v)=|v|^{-B}$ for some $B>0$; then  
the condition (C1) holds for any $B>0$, while  (\ref{covCH}) holds  
if and only if $d\varepsilon +B<4$. Since  for small $\tau >0$, 
$\int_0^{\tau}\rho^{d-1-B-d+4}\, d\rho
=C\, \tau^{4-B}\leq 
I(\tau)\leq \int_0^{\tau}\rho^{d-1-B-d+4-\xi (d-4)^+}\, d\rho
=C\, \tau^{4-B-\xi (d-4)^+}$ 
for any small $\xi >0$, one has $\lim_{\tau\rightarrow 0} I(\tau^{\frac{1}{4}+\nu})^{-1}  
[\tau + \tau^{\frac{1}{2}}\, I(\tau^{\frac{1}{4}-\nu})]=0$ for  
every $B>0$ and $0<\nu <\frac{B\wedge 1}{16}$ , while  
(\ref{Ielip}) holds if and only if $B+{d\varepsilon}  
> \frac{7}{2}\vee \left(4-\left(\frac{{d\varepsilon}}{4}\wedge  
\frac{a}{2}\right)\right)$.  
\end{rk} 
\begin{rk}\label{petitd} The proofs of Theorems 1.5 in \cite{CW1} and
Theorems 1.2-1.4 in \cite{CW2} extend to the case of the space-time
white noise $F$ in dimension $d\leq 3$ under the homogeneous
  Dirichlet boundary conditions on $[0,\pi]^d$.
\end{rk} 
%\bpr 
{\it Proof of Theorem \ref{densiteCH}:} The proof, which is similar to that   
of Theorem 1.2 in  
\cite{CW2} and Theorem 3.1 in \cite{MS} is only sketched in case  
(ii). According to Theorem \ref{holdCH}, the trajectories of the solution $u(n)$ to  
(\ref{u(n)})  almost surely belong to  
$\mathcal{C}^{\lambda,\mu}([\frac{t_0}{2},T]\times Q)$ for  
$0<\lambda<\frac{1}{2}\wedge \frac{d \varepsilon}{8}$ and  
$0<\mu<1\wedge \frac{d\varepsilon}{2}\wedge a$. (Note that  
according to Theorem \ref{holdCH},  the trajectories of $u$ have  
the same H\"older regularity.) Using Theorem 2.1.2 and the  
following remark in \cite{N}, it suffices to show that for every  
$n\geq 1$ and $M\geq 1$, the $l\times l$ Malliavin covariance  
matrix $\Gamma(n)$ defined by  
\[ \Gamma(n)(i,j)=<Du(n)(t_0,x_i),Du_n(t_0,x_j)>_{\mathcal{H}_T}\]  
is almost surely invertible on the set  
$\tilde{\Omega}(M)=\cap_{i=1}^l \{|\sigma(u(n)(t_0,x_i))|\geq  
\frac{1}{M}\}.$ As usual, this reduces to proving that for any  
vector $v\in \mathbb{R}^l$, with $|v|=1$,  
$<\Gamma(n)v,v>_{\mathbb{R}^l}>0$ a.s. on $\tilde{\Omega}(M)$.  
  
For $1\leq i\leq l$, $r\leq t$, (\ref{Dv}) in Proposition  
\ref{Dlip} shows that  
$D_{r,*}u(n)(t_0,x_i)=G(t_0,x_i;r,*)$\linebreak[3]$\times  
\sigma(u(n)(r,*))+U(t_0,r,x_i)$; for a fixed unit vector $v$ of  
$\mathbb{R}^d$, a usual argument shows that given $\tau \in ]  
0,\frac{t_0}{2}]$,  
\begin{equation}\label{decomposition}<\Gamma(n)v,v\geq  
\frac{I_1}{4}+\frac{1}{4}\, \sum_{i=1}^l\sum_{j\neq i; j=1}^l  
v_i\, v_j\, I_2(i,j)-\frac{l}{2}\sum_{i=1}^l v_i^2\, I_3(i)-l\,  
\sum_{i=1}^l v_i^2\, I_4(i)\, ,\end{equation} where  
\begin{eqnarray*}  
I_1&=&\sum_{i=1}^l v_i^2\, \int_{t_0-\tau}^{t_0}\|G(t_0,x_i;r,y)\,  
\sigma(u(n)(r,x_i))\|_{\mathcal{H}}^2\, dr\, ,\\  
I_2(i,j)&=&\int_{t_0-\tau}^{t_0} <G(t_0,x_i;r,*)\,  
\sigma(u(n)(r,x_i)),G(t_0,x_j;r,*)\,  
\sigma(u(n)(r,x_j))>_{\mathcal{H}}\, dr\, ,\\  
I_3(i)&=&\int_{t_0-\tau}^{t_0} \Big\|G(t_0,x_i;r,y)\, \Big[  
\sigma(u(n)(r,y))- \sigma(u(n)(r,x_i))\Big\|_{\mathcal{H}}^2\,  
dr\, ,\\  
I_4(i)&=&\int_{t_0-\tau}^{t_0}\|U(t_0,r,x_i)\|_{\mathcal{H}}^2\,  
dr\, .  
\end{eqnarray*} 
Remark \ref{petrov} shows that condition (\ref{minodiag}) is satisfied.  
Let $\bar{c}=\inf\{d(x_i,\partial Q)\, ,\, 1\leq i\leq l\}$ and  
suppose that for the constant $C_2$ defined in the proof of Lemma  
\ref{CNS} (ii), $2\, C_2\, \tau^{\frac{1}{4}}\leq \bar{c}$. Then  
(\ref{minoCNS=}) and (\ref{minoCNSdiff}) imply that on  
$\tilde{\Omega}(M)$, for $\nu >0$, $\tau$ small enough and  
$d=4,5$,  
\begin{equation}\label{minoI1} I_1\geq C\, (\sum_{i=1}^l |v_i|^2)\,  
\frac{1}{M^2}\, \int_{B_d\left(0, C_2\,  
\tau^{\frac{1}{4}+\nu}\right)} f(v)\, |v|^{-d+4}\,  
\ln\big(|v|^{-1}\big)^{(5-d)^+}\, dv\, .  
\end{equation}  
We now prove upper estimates of $I_2(i,j)$ up to $I_4(i)$. Let  
$m=\inf\{|x_i-x_j|\, ,\, 1\leq i<j\leq l\}$, let $c_1$ and $C_1$  
denote the constants appearing in condition (C1), and let $k\in  
]0,\frac{1}{3}[$ be such that $k\, (1+c_1)<\frac{1}{2}$. Fix  
$1\leq i<j\leq l$; then $I_2(i,j)\leq C\, \|\sigma\|_{\infty}^2\,  
J_2(i,j)$, where  
\[J_2(i,j)=\int_0^\tau r^{-\frac{d}{2}}\, dr \int_Q dy\int_Q  
dz \exp\left(-c  
\frac{|x_i-y|^{\frac{4}{3}}}{r^{\frac{1}{3}}}\right)\,  
f(y-z)\,\exp\left(  
-c\frac{|x_j-z|^{\frac{4}{3}}}{r^{\frac{1}{3}}}\right)\, .\] We  
split the integral on $Q\times Q$ in several parts. Indeed, if  
$|y-z|\geq k m$, the continuity of $f$ implies that $f(y-z)\leq  
C<+\infty$. Suppose now that, $\vert y-z\vert\leq km$. Then  if $|y-x_i|\leq c_1\, |y-z|$,  (C1)  
implies that $f(y-z)\leq C_1\, f(y-x_i)$, while  
$|z-x_j|\geq\Big|\, |x_i-x_j|-(|x_i-y|+|y-z|)\, \Big|\geq m\,  
\Big(1-k(1+c_1)\Big)\geq \frac{m}{2}$. Similarly, if $|z-x_j|\leq  
c_1\, |y-z|$, then $f(y-z)\leq C_1 f(z-x_j)$ and $|y-x_i|\geq  
\frac{m}{2}$. Finally, suppose that $|y-z|\leq k m$, $c_1\,  
|y-z|\leq |y-x_i|\wedge |z-x_j|$; then since $k<\frac{1}{3}$, one  
of the norms $|y-x_i|$ or $|z-x_j|$ (say $|z-x_j|$) is larger than  
$\frac{m}{3}$. Thus, $J_2(i,j)\leq J_{2,1}(i,j)+2\,  
J_{2,2}(i,j)+2\, J_{2,3}(i,j)$, where  
\begin{eqnarray*}  
J_{2,1}(i,j)&=&\int_0^\tau r^{-\frac{d}{2}}\int\int_{|y-z|\geq  
km}\exp\left(-c\frac{|y-x_i|^{\frac{4}{3}}}{r^{\frac{1}{3}}}\right)\,  
f(y-z)\,\exp\left(-c\frac{|z-x_j|^{\frac{4}{3}}}{r^{\frac{1}{3}}}\right)\,  
dydz\leq C\, \tau\, ,  
\\ J_{2,2}(i,j)&\leq & C_1\,  
\int_0^\tau r^{-\frac{d}{4}}\, dr \int_{B_d(0,R)}f(v)\,  
\exp\left(-c\frac{|v|^{\frac{4}{3}}}{r^{\frac{1}{3}}}\right)\, dv  
\int_{|z|\geq \frac{m}{2} r^{-\frac{1}{4}}} \exp\left(  
-c|z|^{\frac{4}{3}}\right) \, dz\\ &\leq& C_1\int_{B_d(0,R)}  
f(v)\, dv\int_0^\tau  
r^{-\frac{d}{4}}\,\exp\left(-c\frac{|v|^{\frac{4}{3}}}{r^{\frac{1}{3}}}  
\right)\, \exp\left( -\tilde{c} r^{-\frac{1}{3}}\right)\,  dr\\&\leq  
&C \, \exp(-\tilde{c}\tau^{-\frac{1}{3}})\, \int_{B_d(0,R)} f(v)\,  
|v|^{-d+4}\,\ln\big(|v|^{-1}\big)^{(5-d)^+}\, dv\leq C \,  
\exp(-\tilde{c}\tau^{-\frac{1}{3}})\, ,\\ J_{2,3}(i,j)&\leq&C_1\,  
\int_0^\tau r^{-\frac{d}{4}}\, dr\int_Q dy\int_{|z-x_j|\geq  
\frac{m}{3}}  f(y-z)  
\exp\left(-\tilde{c}\frac{|y-z|^{\frac{4}{3}}}{r^{\frac{1}{3}}}\right)  
\exp\left(-\tilde{c}\frac{|z-x_j|^{\frac{4}{3}}}{r^{\frac{1}{3}}}\right)\,  
dz\\&\leq& C \,\exp(-\tilde{c}\tau^{-\frac{1}{3}})\,  .  
\end{eqnarray*}  
Hence, for $\tau$ small enough,  
\begin{equation}\label{majoI2} \Big|\sum_{i=1}^l\sum_{j\neq i; j=1}^l  
v_i\, v_j\, I_2(i,j)\Big|\leq C\, \tau\, .  
\end{equation}  
Fubini's theorem, the Lipschitz property of $\sigma$, the H\"older  
regularity of $x\rightarrow u(n)(t,x)$ uniformly for  
$\frac{t_0}{2}\leq t\leq t_0$, Schwarz's inequality and  
$|y-z|^{\frac{4}{3}}\leq 2^{\frac{1}{3}}\Big(  
|y-x_i|^{\frac{4}{3}} +|z-x_j|^{\frac{4}{3}}\Big)$ yield for any  
$i\leq l$, $0<\tau\leq \frac{t_0}{2}$:  
\begin{eqnarray*}  
E(|I_3(i)|)&\leq& C\int_{t_0-\tau}^{t_0} \int_Q \int_Q dy \, dz\,  
|G(t,x_i;r,y)|\, f(y-z)\, |G(t,x_i;r,z)|\\&&\qquad \times  
 E\Big(|u(n)(r,y)-u(n)(r,x_i)| |u(n)(r,z)-u(n)(r,x_i)|\Big)dr \\ &\leq&C\,  
\int_0^{\tau} r^{-\frac{d}{2}} dr\int_Q \int_Q dy dz  
\exp\left(-\tilde{c}\frac{|y-z|^{\frac{4}{3}}}{r^{\frac{1}{3}}}\right)\,  
f(y-z)\, |z-x_i|^{\frac{\mu}{2}} \exp\left(  
-\tilde{c}\frac{|z-x_i|^{\frac{4}{3}}}{r^{\frac{1}{3}}}\right)  
\\&\leq & C \int_{B_d(0,R)} f(v)\, dv \int_0^{\tau}  
r^{-\frac{d}{4}+\frac{\mu}{8}}\, \exp\left(  
-\tilde{c}\frac{|v|^{\frac{4}{3}}}{r^{\frac{1}{3}}}\right) \,  
dr\int_{B_d(0,R)} |z|^{\frac{\mu}{2}}  
\exp(-\tilde{c}|z|^{\frac{4}{3}})\, dz  
\\ &\leq &C \int_{B_d(0,R)}f(v)\, |v|^{-d+4+\frac{\mu}{2}}\, dv  
  \int_{|v|^{\frac{4}{3}}\,\tau^{-\frac{1}{3}}}^{+\infty}\,  
 s^{-4+\frac{3d}{4}-\frac{3\mu}{8}}\,  
\exp(-\tilde{c} s)\, ds\, ,  
\end{eqnarray*}  
where to obtain the last integral, we have set  
$s=|v|^{\frac{4}{3}}\, r^{-\frac{1}{3}}$. We split the last  
integral on $\{|v|\leq \tau^{\frac{1}{4}-\nu}\}$ and its  
complement for some $\nu >0$. Then using (\ref{covCH}), a  
straightforward computation yields for $d=4,5$ and $\tau$ small  
enough:  
\begin{eqnarray*}  
E(|I_3(i)|)&\leq& C\int_{B_d(0,\tau^{\frac{1}{4}-\nu})} f(v)  
|v|^{4-d+\frac{\mu}{2}}\, \Big(1+|v|^{-\frac{\mu}{2}}\,  
\tau^{\frac{\mu}{8}}\Big)\, dv\\ && +C \int_{B_d(0,R)}\,  
1_{\{|v|\geq \tau^{\frac{1}{4}-\nu}\}}\,  f(v)  
|v|^{4-d+\frac{\mu}{2}}\,\exp(-\tilde{c}\tau^{-\frac{4\nu}{3}} )\,  
dv\\ &\leq &C\, \tau^{\frac{d\varepsilon}{4}+\frac{\mu}{8}-\nu  
(d\varepsilon +\frac{\mu}{2})}\, .  
\end{eqnarray*} Thus, choosing $\nu$ close enough to 0 and  
$\tau_0$ small enough, we deduce that for $0<\tau \leq \tau_0$:  
\begin{equation}\label{majoI3} E(|I_3(i)|)\leq C \tau^A\; ,  
\quad 0<A<\frac{d\varepsilon}{4}+\frac{\mu}{8}\, .  
\end{equation}  
Finally, using the decomposition of $U(t_0,r,x_i)$, and Fubini's  
theorem, we obtain that $E(|I_4(i)|)\leq \sum_{i=1}^4 T_j$, where  
\begin{eqnarray*}  
T_1&=&\int_{t_0-\tau}^{t_0} dr E\left( \left\| \int_r^{t_0}\int_Q  
G(s,x_i;r,y)\, S_i(n)(s,y)\, D_{r,*}u(n)(s,y)\, F(ds,dy)  
\right\|^2_{\mathcal{H}}\right) \, ,\\  
 T_2&=&\int_{t_0-\tau}^{t_0} dr E\left( \left\|  
\int_r^{t_0}\int_Q \Delta G(s,x_i;r,y)\, R'_n(u(n)(s,y))\,  
D_{r,*}u(n)(s,y)\, dy ds \right\|^2_{\mathcal{H}}\right) \, ,\\  
 T_3&=&\int_{t_0-\tau}^{t_0} dr E\left( \left\|  
\int_r^{t_0}\int_Q  G(s,x_i;r,y)\, \partial_3 g(s,y,(u(n)(s,y))\,  
D_{r,*}u(n)(s,y)\, dy ds \right\|^2_{\mathcal{H}}\right) \, ,\\  
 T_4&=&\sum_{j=1}^N \int_{t_0-\tau}^{t_0} dr E\left( \left\|  
\int_r^{t_0}\int_Q H_j(s,x_i;r,y)\, \partial_3  
b_j(s,y,(u(n)(s,y))\, D_{r,*}u(n)(s,y)\, dy ds  
\right\|^2_{\mathcal{H}}\right) \,.  
\end{eqnarray*}  
The isometry property for Hilbert-space valued martingales,  
(\ref{majDv}) with $p=1$, Schwarz's inequality and  
(\ref{intpsiTdiff}) or (\ref{intpsiT=}) imply for $\nu$ small enough  
\begin{eqnarray}\nonumber T_1&\leq &C \int_{t_0-\tau}^{t_0}  
ds \int_Q \int_Q dy dz \,  |G(s,x_i;r,y)|\, f(y-z)\,  
|G(s,x_i;r,z)|\\\nonumber &&\qquad \times  \sup_{(s,y)\in  
[t_0-\tau , t_0]\times Q} E\left(\int_{t_0-\tau}^s  
\|D_{r,*}u(n)(s,y)\|^2_{\mathcal{H}}  dr\right)\\&\leq& C\, \tau  
\left[ I(\tau^{\frac{1}{4}-\nu})+\exp\left(-\tilde{c}  
\tau^{\frac{3\nu}{4}}\right)\right]\, . \label{majoT1}  
\end{eqnarray}  
Minkowski's and Schwarz's inequalities, then Fubini's theorem and  
(\ref{majoDvst}) imply  
\begin{eqnarray}\nonumber  
T_2&\leq& C\int_{t_0-\tau}^{t_0} (t_0-r)^{\frac{1}{2}}\, dr  
\int_r^{t_0}\int_Q \Delta G(t_0,x_i;s,y)\, R'_n(u(n)(s,y))\,  
E(\|D_{r,*} u(n)(s,y))\|_{\mathcal{H}}^2)\, dy ds\\\nonumber  
 &\leq& C_n \,  
\tau^{\frac{1}{2}}\int_{t_0-\tau}^{t_0} ds \int_Q dy \, \Delta  
G(t_0,x_i;s,y)\, \sup_{(s,y)\in [t_0-\tau , t_0]\times Q}\,  
E\left( \int_{t_0-\tau}^s \, \|D_{r,*}  
u(n)(s,y))\|_{\mathcal{H}}^2\, dr\right)\\&\leq& C(n)\, \tau  
\left[ I(\tau)+  
  \tau^{\frac{1}{2}}\right] \, . \label{majoT2}  
\end{eqnarray}  
A similar computation, using Minkowski's inequality, Fubini's  
theorem, Lemma \ref{convolution} with $\rho=\infty$ and $q=1$ and  
(\ref{majoDvst}) yields  
\begin{eqnarray}\nonumber  
 T_3&\leq& C_n\, \tau\,  
\int_{t_0-\tau}^{t_0} ds \int_Q dy \, G(t_0,x_i;s,y)\,  
\sup_{(s,y)\in [t_0-\tau , t_0]\times Q}\,  E\left(  
\int_{t_0-\tau}^s \, \|D_{r,*} u(n)(s,y))\|_{\mathcal{H}}^2\,  
dr\right)\\ \label{majoT3}&\leq& C_n\, \tau^{2}\, \left[  
I(\tau)+  
\tau^{\frac{1}{2}}\right]\, ,\\  
 \nonumber  T_4&\leq& C_n  
\tau^{\frac{1}{4}} \sum_{j=1}^N \int_{t_0-\tau}^{t_0} ds \int_Q dy  
H_j(t_0,x_i;s,y) \sup_{(s,y)\in [t_0-\tau , t_0]\times Q}\,  
E\left( \int_{t_0-\tau}^s  \|D_{r,*} u(n)(s,y))\|_{\mathcal{H}}^2  
dr\right)\\  
 \label{majoT4}  &\leq& C_n\, \tau^{\frac{1}{2}}\,  
 \left[ I(\tau)+  
\tau^{\frac{1}{2}}\right]\, ,  
\end{eqnarray}  
The inequalities (\ref{majoT1})-(\ref{majoT4}) yield that  
\begin{equation}\label{majoI4}  
E(|I_4(i)|)\leq C(n)\,  \left[\tau^{\frac{1}{2}}\,  
  I(\tau^{\frac{1}{4}-\nu})+  
\tau\right]\, .  
\end{equation}  
Finally, the inequalities (\ref{decomposition})-(\ref{majoI3}) and  
(\ref{majoI4}) imply that for $\rho >0$ such that $\tau+\tau^\rho \leq  
\frac{1}{2} I(\tau^{\frac{1}{4}+\nu})$ (which exists because of  
(\ref{Ielip})), we have for small enough $\tau$:  
\begin{eqnarray*}  
P(<v,\Gamma(n)v>_{\mathbb{R}^l}<\tau^\rho)&\leq&P\left(  
I_3+I_4\geq C\,  [ I( \tau^{\frac{1}{4}+\nu}) -\tau  
-\tau^\rho]\right)\\&\leq &P\left(I_3+I_4\geq \frac{1}{2} \, I(  
\tau^{\frac{1}{4}+\nu})\right) \leq C\,  I(  
\tau^{\frac{1}{4}+\nu})^{-1} \left[ \tau^A + \tau^{\frac{1}{2}} \,  
I(\tau^{\frac{1}{4}-\nu})\right]  
\end{eqnarray*}  
for some $A<\frac{\mu}{8}+\frac{d\varepsilon}{2}$; then the  
condition (\ref{Ielip}) concludes the proof of (ii). \epr  
  
%\newpage 

\end{document}